\documentclass{article}
\usepackage[utf8]{inputenc}
\usepackage{amsmath}
\usepackage{amsthm}
\usepackage{amsfonts}
\usepackage{amssymb}
\usepackage{hyperref}
\usepackage{listings}
\usepackage{xcolor}

\hypersetup{
    colorlinks,
    linkcolor={red!50!black},
    citecolor={blue!50!black},
    urlcolor={blue!20!black}
}

\usepackage{float}
\usepackage{geometry}
\usepackage{graphicx}
\usepackage{subcaption}
\usepackage{caption}
\usepackage{xspace}

\newcommand{\treeA}{{\rt[[],[[],[[]]]]}}
\newcommand{\treeB}{{\rt[[[]],[[],[]]]}}
\newcommand{\treeC}{{\rt[[[],[[],[]]]]}}

\newcommand{\order}{{\mathcal O}}
\newcommand{\coeff}{u}
\newcommand{\coeffb}{v}
\newcommand{\spine}{trunk\xspace}
\newcommand{\remainder}{remainder\xspace}

\newcommand{\tree}{\tau}
\newcommand{\stepsizex}{h_0}

\usepackage{epsfig}
\usepackage{forest}

\newdimen\myshiftaligndimen

\forestset{
  declare dimen register={shift align},
  shift align'=0pt,
  default preamble={
    /tikz/baseline={([yshift=-2.5pt]current bounding box.center)},
    for tree={%
      grow=north,
      circle,
      fill,
      minimum size=2.5pt,
      inner sep=0pt,
    },
  },
  */.style={
    delay+={append={[]},}
  },
  rooted tree/.style={
    for tree={
      grow'=90,
      s sep=2.5pt,
      delay={
        if content={*}{
          content=,
          append={[]}
        }{}
      }
    },
    before typesetting nodes={
      for tree={
        circle,
        fill,
        minimum width=3pt,
        inner sep=0pt,
        child anchor=center,
      },
    },
    before computing xy={
      for tree={
        l=5pt,
      }
    }
  },
  big rooted tree/.style={
    for tree={
      grow'=90,
      parent anchor=center,
      child anchor=center,
      s sep=2.5pt,
      if level=0{
        baseline
      }{},
      delay={
        if content={*}{
          content=,
          append={[]}
        }{}
      }
    },
    edge={
        semithick,
        -Latex
    },
    before typesetting nodes={
      for tree={
        circle,
        fill,
        minimum width=5pt,
        inner sep=0pt,
        child anchor=center,
      },
    },
    before computing xy={
      for tree={
        l=15pt,
      }
    }
  },
  my node/.style={
    circle, fill=#1,
  },
  gray node/.style={
    my node=gray, inner sep=#1mm,
  },
  red node/.style={
    my node=red, inner sep=#1mm,
  },
  blue node/.style={
    my node=blue, inner sep=#1mm,
  },
  green node/.style={
    my node=green, inner sep=#1mm,
  },
}
\DeclareDocumentCommand\rt{o}{\Forest{rooted tree [#1]}}
\DeclareDocumentCommand\bigrt{o}{\Forest{big rooted tree [#1]}}

\newtheorem{thm}{Theorem}

\newtheorem{example}{Example}
\newtheorem{remark}{Remark}
\newtheorem{dfn}{Definition}

\title{Pseudo-energy-preserving explicit Runge-Kutta methods}
\author{Gabriel A. Barrios de Le\'on\thanks{Escuela de Ciencias F\'isicas y Matem\'aticas, Universidad de San Carlos de Guatemala} \and David I. Ketcheson\thanks{Applied Mathematics and Computational Science, CEMSE Division, King Abdullah University of Science and Technology (KAUST), Thuwal, 23955-6900, Kingdom of Saudi Arabia} \and Hendrik Ranocha\thanks{Institute of Mathematics, Johannes Gutenberg University Mainz, Staudingerweg 9 55130 Mainz, Germany}}

\begin{document}
\date{}
\maketitle

\begin{abstract}
Using a recent characterization of energy-preserving B-series, we derive the
explicit conditions on the coefficients of a Runge-Kutta method that ensure
energy preservation (for Hamiltonian systems) up to a given order in the step
size, which we refer to as the pseudo-energy-preserving (PEP) order.  We study
explicit Runge-Kutta methods with PEP order higher than their classical order.
We provide examples of such methods up to PEP order six,
and test them on Hamiltonian ODE and PDE systems.  We find that these
methods behave similarly to exactly energy-conservative methods over
moderate time intervals and exhibit significantly smaller errors, relative to
other Runge-Kutta methods of the same order, for moderately long-time simulations.
\end{abstract}

\section{Introduction}
The preservation of energy is an essential concern in the numerical integration
of Hamiltonian systems.  Decomposing the numerical
error into a component along the energy-preserving manifold and a component that moves
off of the manifold, it is generally accepted that the non-energy-preserving component
of the error is more harmful, especially in long-term simulations.  This component both
allows for qualitatively wrong behavior and leads to compounding errors that result
in greater overall quantitative error.
It is therefore valuable to develop and use methods that exactly preserve energy.

Traditional Runge-Kutta methods cannot be exactly energy-preserving
for all Hamiltonian systems \cite{celledoni2009energy}, although we note that
even explicit methods can be exactly energy-preserving for a specific
Hamiltonian system (see e.g. \cite[Section~5.1]{ranocha2020energy}.)
Exact energy preservation for arbitrary Hamiltonian systems can be achieved with (implicit) continuous-stage
Runge-Kutta methods \cite{Hairer2010}, or by other approaches such as projection
\cite[Section~IV.4]{hairer2006geometric} or
relaxation \cite{ketcheson2019relaxation,ranocha2020relaxation,ranocha2020general}.

In this work we explore an alternative, in which the order of the
energy conservation error is made as high as possible.  While this
does not lead to exact energy conservation, it can ensure that the
error in the energy is much smaller than the overall error of the method,
and is acceptably small for relatively long times.  Even exactly
energy-conserving methods will violate energy conservation when
implemented in floating-point arithmetic, so numerical energy conservation
is effectively a matter of making the energy error small in any case.

In this work we develop explicit Runge-Kutta methods with maximal
\emph{pseudo-energy-preserving} (PEP) order.
To understand what is meant by PEP order, consider a method whose
classical order of accuracy is $p$; i.e.
the error of the numerical solution in one step is $\order(h^{p+1})$. Then
the error of the Hamiltonian after one step is also $\order(h^{p+1})$ in
general. We say that the method has PEP order $q$ if the error of the
Hamiltonian in one step is $\order(h^{q+1})$, with $q  \ge p$.

We focus on explicit Runge-Kutta methods.
To analyze the PEP order of a method, we consider the B-series of the modified equation of a
Runge-Kutta method. Theorems~2--3 of \cite{celledoni2010energy} provide a means
to write a set of necessary and sufficient conditions on the coefficients of a B-series in
order for it to be energy-preserving. A method then has PEP order $q$ if these
algebraic conditions are satisfied for the coefficients of all terms up to (and including)
order $q$.

Using the relationship between the modified equation and the B-series of the method itself,
and then the relationship of this series with the Butcher coefficients, one can obtain
energy-preservation conditions directly in terms of the coefficients of the RK method.
We work out these conditions in Section \ref{sec:conditions}.  In Section \ref{sec:methods}
we solve these conditions for explicit RK methods with a given number of stages and prescribed
classical order of accuracy, maximizing the PEP order subject to these constraints.  In Section
\ref{sec:tests}, we apply these methods to a range of Hamiltonian systems.  We study in
particular the long-time behavior of the global error.

\subsection{Related ideas}
Structure-preserving B-series are known to improve numerical simulations by maintaining certain invariants exactly.
Examples of these well-studied structures are volume, linear quadratic invariants, or the symplectic form \cite{munthe2016aromatic,hong2014invariants,bogfjellmo2019algebraic}.
Pseudo-structure-preserving methods, in contrast, aim to balance efficiency with approximate preservation depending on the problem.
In the context of stochastic partial differential equations (PDEs), the classical concept of energy is replaced by a modified energy
as in \cite{debussche2009modenergy}, and the development of integrators that preserve modified energy
is an active area of research as well.

Numerical energy preservation and symplecticity are closely related, and a number
of previous works have studied so-called pseudo-symplectic Runge-Kutta methods
\cite{aubry1998pseudo,calvo2010approximate,capuano2017explicit,stepanov2023eight}.
These are explicit schemes that are constructed to satisfy the conditions for
symplecticity (equivalent to the preservation of quadratic invariants) up to
some order higher than the classical order of accuracy.
A similar approach has been used to develop the conditions for an aromatic B-series method to be volume-preserving \cite{Laurent2023aromatic} and to devise
pseudo-volume-preserving methods \cite{munthe2016aromatic,bogfjellmo2019algebraic}.
In the numerical experiments in this work we provide some comparisons of pseudo-energy-preserving
and pseudo-symplectic methods in practice.  An interesting area for future work might
be a more complete theoretical and practical comparison of the conditions among these three classes of methods.

In order to preserve energy exactly, the class of continuous-stage Runge-Kutta (CSRK) methods
has been introduced, starting with the average vector field (AVF) method
\cite{mclachlan1999geometric}.  This is a second-order method that preserves energy
for any Hamiltonian system but requires the solution of a system of algebraic equations
similar to what is required for a fully-implicit Runge-Kutta method.
Higher-order accurate energy-preserving CSRK methods have been introduced
by Miyatake \& Butcher \cite{Miyatake2014,miyatake2016characterization}.

Alternatively, by applying a projection after each step, any Runge-Kutta method
can be made to preserve a specific desired functional, such as the energy.
A certain kind of oblique projection (in the direction of the step) has the
advantage of also preserving linear invariants; see
e.g.  \cite{calvo2006preservation,calvo2010approximate,Ranocha2020RRKHamiltonian}.
In a similar vein, the discrete gradients approach can be used to enforce conservation
of a prescribed nonlinear functional \cite{mclachlan1999geometric}.

\subsection{B-series and Runge-Kutta methods}
In this work we make use of B-series \cite{Butcher2003,hairer2006geometric,butcher2021b}, written in the form
\begin{align}
    B(\coeff, hf, y) & := \sum_{\tree \in T} \frac{h^{|\tree|}}{\sigma(\tree)} \coeff(\tree) F_f(\tree)(y),
\end{align}%
where $T$ is the set of all rooted trees, $\sigma(\tree)$ is the \emph{symmetry} of $\tree$,
and $|\tree|$ is the \emph{order} of
the rooted tree $\tree$, i.e., the number of nodes in $\tree$.
The factor $F_f(\tree)(y)$ involves the derivatives of $f$ and is referred to
as an \emph{elementary differential}.

To every Runge-Kutta method we associate two B-series.  First, the map, with
coefficients $\coeff(\tree)$ such that the numerical solution obtained by applying
the method with step size $h$ is
\begin{align} \label{Bmap}
    y_{n+1} & = y_n + B(\coeff, hf, y_n).
\end{align}%
The coefficients $\coeff(\tree)$ are known as the \emph{elementary weights} of the
method, and depend only on the RK coefficients $A, b$.
Secondly we associate to an RK method the flow, with coefficients $v(\tree)$
defined such that the numerical solution obtained by applying the method with
step size $h$ is the exact solution of the ODE
\begin{align} \label{Bflow}
    \dot{y}(t) & = B(\coeffb, hf, y).
\end{align}%
This ODE \eqref{Bflow} is known as the \emph{modified equation}
\cite[Chapter~IX]{hairer2006geometric}.

The coefficients $\coeff$ and $\coeffb$ of the map and flow of a RK method
are related in such a way that either can be determined from the other in a
sequential process starting with the coefficients of lower-order
trees and proceeding to those of higher order \cite{Chartier2010}.
This process results from applying the substitution law for B-series and writing down the resulting equations for each tree.
This is detailed further in Section \ref{sec:map-flow}.

\section{Pseudo-energy-preservation conditions for Runge-Kutta methods applied to Hamiltonian systems} \label{sec:conditions}
Let us consider the application of a Runge-Kutta method to a Hamiltonian system of ODEs:
\begin{align} \label{hamiltonian-system}
    \dot y(t) = J \nabla H(y).
\end{align}
Here $y \in \mathbb{R}^d$, the \emph{Hamiltonian} $H(y) \in C^q(\mathbb{R}^d)$,  and
$J$ is a constant non-singular skew-symmetric matrix.  For every solution of \eqref{hamiltonian-system},
$H(y)$ is constant in time.  We say a Runge-Kutta method is \emph{energy-preserving} if
it inherits this property:

\begin{dfn}
A Runge-Kutta method is energy-preserving
if it yields $H(y_{n+1})=H(y_n)$ for every Hamiltonian system \eqref{hamiltonian-system}.
\end{dfn}

\begin{dfn}
A Runge-Kutta method is \emph{pseudo-energy-preserving} of order $(p,q)$ if it has classical order
$p$ and satisfies $H(y_{n+1})=H(y_n) + \order(h^{q+1})$ for every Hamiltonian system \eqref{hamiltonian-system}.
\end{dfn}

Throughout this work, when discussing methods with pseudo-energy-preserving order $q$,
we assume that $H$ is at least $q$ times differentiable.
Necessary and sufficient conditions for pseudo-energy-preserving methods have been given
in terms of the method's modified equation in \cite[Theorems~2-3]{celledoni2010energy}.
Here we first recall this result and provide an explanation and examples of the explicit conditions for energy preservation.
We then derive the corresponding conditions on the coefficients of a pseudo-energy-preserving
RK method.

\subsection{Energy-preserving B-series flows}
Let $T_p$ denote the set of rooted trees of order at most $p$.
Given a tree $\tree$, let $L(\tree)$ denote the set of its leaves (nodes with no children).  For each leaf
$\ell\in L(\tree)$, define the corresponding \spine $S(\ell)$ as the ordered set of nodes in a direct path starting from the
root and ending at the leaf.  When the \spine (the nodes and adjacent edges) is removed from $\tree$, what remains is a
sequence of forests, which is also ordered according to the node that each forest was
attached to; we will refer to this sequence of forests as the \remainder $\tree \setminus S(\ell) = \{f_1, f_2, \dots, f_m\}$.
Finally, we define the \emph{energy-preserving conjugate} $\tree^*(\ell)$ of a tree $\tree$ and leaf $\ell$ as follows.
We take the \spine and reattach the forests in $\tree \setminus S(\ell)$ but in reverse order; i.e. $f_m$ is attached to the root,
$f_{m-1}$ is attached to the second node of the \spine, etc., until $f_1$ is attached to the second-to-last
node of the \spine.  Examples are given in Tables~\ref{tbl:example}, \ref{tbl:example2},
and \ref{tbl:example3}.

An alternative definition of $\tree^*(\ell)$ is as follows
(see the proof of \cite[Thm.~2]{celledoni2010energy}).  Given the tree
$\tau$ with root $r$ and leaf $\ell$, relabel the parent of $\ell$ as the root.
Then remove $\ell$ and reattach it as a child of the original root node $r$.
The resulting rooted tree is $\tau^*(\ell)$.

\begin{table}[htbp]
    \centering
    \caption{Example of a tree $\tree$ and its energy-preserving conjugates.\label{tbl:example}}
    \begin{tabular}{lcccc} \hline
    $\tree$ & \rt[[.,red node],[[.,blue node],[[.,green node]]]] \\ \hline
    $\ell$  & \rt[.,red node] & \rt[.,blue node] & \rt[., green node] \\ \hline
    $m$ & 1 & 2 & 3 \\ \hline
    $\tree \setminus S(\ell)$ & \{ \rt[[],[[]]] \}   & \{ \rt[], \rt[[]] \} &  $\{ \rt[], \rt[], \varnothing \}$ \\ \hline
    $\tree^*(\ell)$   & \rt[[],[[],[[]]]] & \rt[[[]],[[],[]]] &    \rt[,[[],[[][]]]] \\ \hline
    \end{tabular}
\end{table}

\begin{table}[htbp]
  \centering
  \caption{Example of a tree $\tree$ and its energy-preserving conjugates.\label{tbl:example2}}
  \begin{tabular}{lcccc} \hline
  $\tree$ & \rt[[.,red node] [.,blue node][[.,green node]]] \\ \hline
  $\ell$  & \rt[.,red node] & \rt[.,blue node] & \rt[., green node] \\ \hline
  $m$ & 1 & 1 & 2 \\ \hline
  $\tree \setminus S(\ell)$ & \{ \{ \rt[], \rt[[]] \}  \}   & \{ \{ \rt[], \rt[[]] \} \} &  $\{ \{ \rt[], \rt[] \} , \varnothing \}$ \\ \hline
  $\tree^*(\ell)$   & \rt[[][][[]]] & \rt[[][][[]]] &  \rt[[[][][]]] \\ \hline
  \end{tabular}
\end{table}

\begin{table}[htbp]
  \centering
  \caption{Example of a tree $\tree$ and three of its four energy-preserving conjugates.  The remaining conjugate is $\tau$ itself.
  \label{tbl:example3}}
  \begin{tabular}{lcccc} \hline
  $\tree$ & \rt[[][[., green node]][[[., blue node]][., red node]]] \\ \hline
  $\ell$  & \rt[.,red node] & \rt[.,blue node] & \rt[., green node] \\ \hline
  $m$ & 2 & 3 & 2 \\ \hline
  $\tree \setminus S(\ell)$ & \{ \{ \rt[], \rt[[]] \} , \rt[[]] \}   & $\{ \{ \rt[], \rt[[]] \}, \rt[], \varnothing \} $ &  $\{ \{ \rt[], \rt[[][[]]] \} , \varnothing \}$ \\ \hline
  $\tree^*(\ell)$   & \rt[[[]][[[]][][]]] & \rt[[[][[[]][][]]]] & \rt[[[][][[[]][]]]] \\ \hline
  \end{tabular}
\end{table}

The signed linear combinations of each tree with each of its energy-preserving conjugates
form a basis for the space of energy-preserving B-series.

\begin{thm}[{\cite[Thms.~2--3]{celledoni2010energy}}]
\label{thm:energy-preserving-flow}
    A B-series flow \eqref{Bflow} for the canonical Hamiltonian system \eqref{hamiltonian-system} is pseudo-energy-preserving of order $q$ if and only if there exist coefficients $\mu_{j,\tree}$
    such that
    \begin{align} \label{thm-eqn}
        B(v,hf,y) & = \sum_{\tree\in T_q} h^{|\tree|} \sum_{\ell_j \in L(\tree)} \mu_{j,\tree} \left( F_f(\tree) + (-1)^{m_j(\tree)} F_f(\tree^*(\ell_j)) \right) + \order(h^{q+1}).
    \end{align}
    Here the sequence $\ell_j$ is an ordering of the leaves of $\tree$ and $m_j(\tree)$ is the distance from the
  root to leaf $\ell_j$.
\end{thm}

If the modified equation (flow) corresponding to a RK method is energy preserving to order $q$,
we say the RK method itself is pseudo-energy-preserving of order $q$.

The implications of Theorem \ref{thm:energy-preserving-flow} are particularly
simple for certain classes of trees, as described in the following examples.

\begin{example}[Bushy trees]
All leaves of a bushy tree $\tree \in \{\rt[[]], \rt[[][]], \rt[[][][]], \dots\}$ are
attached directly to the root. Thus, the energy-preserving conjugate $\tree^*(\ell) = \tree$
for all leaves $\ell$.  Since the distance of each leaf from the root is $m = 1$, we have
\begin{align} \label{f-minus-f}
    F_f(\tree) + (-1)^{m_j(\tree)} F_f(\tree^*(\ell_j) = F_f(\tree)- F_f(\tree) & = 0,
\end{align}
so these trees cannot appear in an energy-preserving B-series flow.
\end{example}

\begin{example}[Tall trees]
For any tall tree
\begin{equation*}
  \tree \in \left\{ \rt[], \rt[[]], \rt[[[]]], \dots \right\},
\end{equation*}
there is only one leaf and the \spine is the entire tree.
For tall trees with an even number of nodes, $m$ is odd and \eqref{f-minus-f}
holds, so these trees cannot appear in an energy-preserving B-series flow.
Meanwhile, for tall trees with an odd number of nodes, $m$ is even, so
\begin{align} \label{f-plus-f}
    F_f(\tree) + (-1)^{m_j(\tree)} F_f(\tree^*(\ell_j) = F_f(\tree) + F_f(\tree) & = 2F_f(\tree),
\end{align}
so the coefficients of these trees are allowed to take any value in an energy-preserving
B-series.
\end{example}

\begin{example}[Self-conjugate trees]
There are other trees for which \eqref{f-plus-f} holds (for some leaf $\ell_j$); e.g.,
\begin{equation*}
  \rt[[][[][]]], \quad
  \rt[[[][[][[]]]]].
\end{equation*}
For each of these trees, choosing the rightmost leaf leads to an
energy-preserving conjugate that coincides with the original tree with an
even distance $m$. Thus, their coefficient can be arbitrary in an
energy-preserving B-series flow.
\end{example}

The condition given in Theorem \ref{thm:energy-preserving-flow} is convenient for checking
the energy-preserving order of a given B-series, since the existence of the coefficients
$\mu_{j,\tau}$ can be phrased in terms of a system of linear equations.  However,
it is less convenient for designing energy-preserving methods, since it does not
explicitly give the algebraic conditions that must be satisfied by the B-series coefficients.
The latter can be obtained by considering a linear system of equations for the
unknown coefficients $\mu_{j,\tau}$.  Because the energy-preserving-conjugacy
relation is transitive, the set of trees of a given order can be partitioned into
disjoint subsets, with the trees in each subset being energy-preserving conjugates
of one another.  One algebraic condition is obtained for each of these sets.

\subsection{Conditions on the flow coefficients}\label{sec:flow-cond}

All consistent RK methods have PEP order at least one.  The condition for
PEP order two is
$$
    \coeffb(\rt[[]]) = 0,
$$
which is just the condition for classical order two.  So every
consistent method with $q\ge 2$ also has $p\ge 2$.

For PEP order three we have one additional condition:
$$
    \coeffb(\rt[[],[]]) = 0.
$$%
Thus PEP order 3 is a weaker condition than classical order 3,
as the latter additionally requires $\coeffb\left(\rt[[[]]]\right)=0$.

For PEP order four we have three additional conditions:
\begin{align*}
        \coeffb\left(\rt[[[[]]]]\right) & = 0, &
        \frac{1}{2}\coeffb\left(\rt[[[][]]]\right) & = \coeffb\left(\rt[[],[[]]]\right), &
        \coeffb\left(\rt[[][][]]\right) & = 0.
\end{align*}

The additional conditions for PEP order five are
\begin{align*}
            \frac{1}{2}\coeffb\left(\rt[[[[][]]]]\right) + \coeffb\left(\rt[[[[]]][]]\right) & = 0, &
            \coeffb\left(\rt[[[[]][]]]\right) - \frac{1}{2}\coeffb\left(\rt[[[]][[]]]\right) & = 0, \\
            \coeffb\left(\rt[[[]][][]]\right) - \frac{1}{3}\coeffb\left(\rt[[[][][]]]\right) & = 0, &
            \coeffb\left(\rt[[][][][]]\right) & = 0.
\end{align*}

Up to fifth order, each tree is conjugate to only itself and/or one
other tree.  Beginning at sixth order, there are triplets of trees among which
each is conjugate to the other two.  The conjugacy relations among such
triplets always result in a single constraint for energy-preserving B-series.
As an example, the conditions involving one such triplet are
\begin{multline}
    \coeffb\left(\treeA\right) F_f\mathopen{}\left(\treeA\mathclose{}\right) + \frac{1}{2}\coeffb\left(\treeB\right) F_f\mathopen{}\left(\treeB\mathclose{}\right) + \frac{1}{2}\coeffb\left(\treeC\right) F_f\mathopen{}\left(\treeC\mathclose{}\right) \\
    =
    \mu_1\left(F_f\mathopen{}\left(\treeA\mathclose{}\right) + F_f\mathopen{}\left(\treeB\mathclose{}\right)\right) + \mu_2\left(F_f\mathopen{}\left(\treeB\mathclose{}\right) + F_f\mathopen{}\left(\treeC\mathclose{}\right)\right) + \mu_3\left(F_f\mathopen{}\left(\treeA\mathclose{}\right) - F_f\mathopen{}\left(\treeC\mathclose{}\right)\right).
\end{multline}%
Equating coefficients of each elementary differential yields a system of three linear
equations, but these equations are linearly dependent, and they have a solution only if
$$
\coeffb\left(\treeA \right) + \frac{1}{2}\coeffb\left(\treeC \right) = \frac{1}{2} \coeffb\left(\treeB\right) .
$$%
The PEP order six conditions are then:
\begin{align*}
            \coeffb\left(\rt[[[[[[]]]]]]\right) = \coeffb\left(\rt[[],[[[],[]]]]\right) =
            \coeffb\left(\rt[[[[]],[[]]]]\right) &= \coeffb\left(\rt[[],[],[],[],[]]\right)  = 0, \\
            \frac{1}{2}\coeffb\left(\rt[[[[[],[]]]]]\right) - \coeffb\left(\rt[[],[[[[]]]]]\right) & = 0, &
            \frac{1}{6}\coeffb\left(\rt[[[[],[],[]]]]\right) + \frac{1}{2}\coeffb\left(\rt[[],[],[[[]]]]\right) & = 0, \\
            \frac{1}{2}\coeffb\left(\rt[[[],[],[[]]]]\right) - \frac{1}{2}\coeffb\left(\rt[[],[[]],[[]]]\right) & = 0, &
            \frac{1}{24}\coeffb\left(\rt[[[],[],[],[],[]]]\right) - \frac{1}{6}\coeffb\left(\rt[[],[],[],[[]]]\right) & = 0 ,\\
            \frac{1}{6}\coeffb\left(\rt[[],[[],[],[]]]\right) - \frac{1}{4}\coeffb\left(\rt[[],[],[[],[]]]\right) & = 0, &
            \coeffb\left(\treeA \right)+ \frac{1}{2}\coeffb\left(\treeC\right) - \frac{1}{2} \coeffb\left( \treeB \right) & = 0, \\
            \coeffb\left(\rt[[[[],[[]]]]]\right) - \coeffb\left(\rt[[[],[[[]]]]]\right) + \coeffb\left(\rt[[[]],[[[]]]]\right) & = 0.
\end{align*}

In principle it is possible to work out the energy-preserving conditions for
trees of any order.  We do not list further conditions here, although the code
used to determine the conditions above can easily be used to generate higher-order
conditions \cite{barrios2024pseudoRepro}.

\subsection{Relation between elementary weights of the map and flow}\label{sec:map-flow}
    As mentioned already, the relationship between the coefficients $v$ of the
    modified equation (flow) and the coefficients $u$ of the RK map is given by
    the substitution law; for details we refer to e.g. \cite{Chartier2010}.  In this section
    we simply apply that law.

    Since the tree of order 2 is not allowed to appear in an energy-preserving
    B-series flow, we see that any numerical method that is at least first-order accurate
    and pseudo-energy-preserving of order 2 must in fact be at least second-order
    accurate.  Thus we simplify what follows by considering only B-series of
    methods that have order of accuracy 2 or greater.

    The following code returns the coefficients $v$ of the
    modified equation (flow) in terms of the coefficients $u$ of the RK map, up to order $p$,
    using the Julia package BSeries.jl \cite{ketcheson2023computing}
    and SymPy \cite{meurer2017sympy}:

\begin{lstlisting}
using BSeries, SymPy
order = p # choose an integer for p
series = bseries(order) do t, series
    return symbols("u_$(butcher_representation(t))", real = true)
end
modified_equation(series)
\end{lstlisting}

By solving the resulting system of equations for the coefficients $u$ of the
RK map, we can express them in terms of the coefficients $v$.
To save space, we write these conditions only up to fourth order
(and assuming a method with classical order at least two):

\begin{equation} \label{uv}
    \begin{aligned}
        \coeff\left(\emptyset\right)  &= 1, \quad & \coeff\left(\rt[]\right)  &= 1,\\
        \coeff\left(\rt[[]]\right) & = \frac{1}{2}, \quad &
        \coeff\left(\rt[[[]]]\right) &= \coeffb\left(\rt[[[]]]\right) + \frac{1}{6},\\
        \coeff\left(\rt[[][]]\right)  &= \coeffb\left(\rt[[][]]\right) + \frac{1}{3}, \quad &
        \coeff\left(\rt[[[[]]]]\right)  &= \coeffb\left(\rt[[[[]]]]\right) + \coeffb\left(\rt[[[]]]\right) +\frac{1}{24}, \\
        \coeff\left(\rt[[[][]]]\right)  &= \coeffb\left(\rt[[[][]]]\right) + \frac{1}{2}\coeffb\left(\rt[[][]]\right) + \coeffb\left(\rt[[[]]]\right) + \frac{7}{12}, \\
        \coeff\left(\rt[[[]][]]\right)  &= \coeffb\left(\rt[[[]][]]\right) + \frac{1}{2} \coeffb\left(\rt[[][]]\right) + \frac{1}{2} \coeffb\left(\rt[[[]]]\right) + \frac{1}{6}, \\
        \coeff\left(\rt[[][][]]\right)  &= \coeffb\left(\rt[[][][]]\right) + \frac{3}{2}\coeffb\left(\rt[[][]]\right) + \frac{1}{4}.\\
    \end{aligned}
\end{equation}

    \subsection{Energy-preserving conditions on the RK coefficients}

    Applying the results of the previous two sections,
    we can obtain conditions for energy
    preservation up to a given order directly in terms of the RK coefficients.
    In order to make the equations easier to read, we write them in terms of
    the elementary weights (i.e. the coefficients $\coeff\left(\tau\right)$).  The formulas
    for the weights in terms of the RK coefficients are available in standard
    references such as \cite{Butcher2003}.

    Combining the conditions from Section \ref{sec:flow-cond} with the relations \eqref{uv}, we obtain
    the following conditions.

    PEP order three:
    \begin{equation}
        \coeff\left(\rt[[][]] \right)= \frac{1}{3}.
    \end{equation}

    PEP order four:
    \begin{equation}
        \begin{aligned}
            \coeff\left(\rt[[[[]]]] \right) & =\coeff\left(\rt[[[]]]\right) -\frac{1}{8}, \quad &
            \coeff\left(\rt[[[]][]]\right) -\frac{1}{2}\coeff\left(\rt[[[][]]]\right) & = \frac{1}{12}, \quad &
            \coeff\left(\rt[[][][]]\right) & = \frac{1}{4}. \\
        \end{aligned}
        \label{cond:p2q4}
    \end{equation}

    PEP order five:
    \begin{equation}
        \begin{aligned}
            \coeff\left(\rt[[[[]]][]]\right) + \frac{1}{2}\coeff\left(\rt[[[[][]]]]\right) & = \coeff\left(\rt[[[[]]]]\right)+\frac{1}{2}\coeff\left(\rt[[[][]]]\right) - \frac{1}{2}\coeff\left(\rt[[[]]]\right) +\frac{1}{24}, & \\
            2\coeff\left(\rt[[[[]][]]]\right) -\coeff\left(\rt[[[]][[]]]\right) & = \coeff\left(\rt[[[[]]]]\right) +\coeff\left(\rt[[[][]]]\right) - \coeff\left(\rt[[[]]]\right)+\frac{1}{24} , & \\
            \coeff\left(\rt[[[]][][]]\right) - \frac{1}{3}\coeff\left(\rt[[[][][]]]\right) & = \frac{1}{12}, \quad &
            \coeff\left(\rt[[][][][]]\right) & = \frac{1}{5}. \\
        \end{aligned}
        \label{cond:p2q5}
    \end{equation}

    PEP order six:
\begin{equation}
  \begingroup
  \small
  \begin{aligned}
            &6 \left(3 \coeff\left(\rt[[[]]]\right)+ \coeff\left(\rt[[[]][]]\right)+ 3 \left(\coeff\left(\rt[[[[][]]]]\right)+ 4 \coeff\left(\rt[[[[[]]]][]]\right) \right) \right) = 2 + 3 \coeff\left(\rt[[[][]]]\right)+ 24 \coeff\left(\rt[[[[]]]]\right)+ 36 \coeff\left(\rt[[[[]]][]]\right)+ 36 \coeff\left(\rt[[[[[][]]]]]\right), \\
            &1 - 12 \coeff\left(\rt[[[]]]\right)+ 24 \coeff\left(\rt[[[]][]]\right)- 36 \coeff\left(\rt[[[]][][]]\right)+ 24 \coeff\left(\rt[[[][]]]\right)- 12 \coeff\left(\rt[[[][][]]]\right)+  24 \coeff\left(\rt[[[[]]]]\right)- 72 \coeff\left(\rt[[[[]]][]]\right)+ 72 \coeff\left(\rt[[[[]]][[]]]\right)- 36 \coeff\left(\rt[[[[][]]]]\right)+ 24 \coeff\left(\rt[[[[][][]]]]\right)= 0, \\
            &12 \left(2 \coeff\left(\rt[[[]][]]\right)+ 6 \coeff\left(\rt[[[]][][]]\right)+ 6 \coeff\left(\rt[[[]][[]]]\right)+ 5 \coeff\left(\rt[[[][]]]\right)+ 2 \coeff\left(\rt[[[[]]]]\right)+ 12 \coeff\left(\rt[[[[]][][]]]\right)\right) = 5 + 24 \coeff\left(\rt[[[]]]\right)+ 144 \coeff\left(\rt[[[]][[]][]]\right)+ 72 \coeff\left(\rt[[[][][]]]\right)+ 144 \coeff\left(\rt[[[[]][]]]\right), \\
            &1 + 60 \coeff\left(\rt[[[]][]]\right)+ 120 \coeff\left(\rt[[[]][][][]]\right)+ 60 \coeff\left(\rt[[[][][]]]\right)=  30 \left(6 \coeff\left(\rt[[[]][][]]\right)+ \coeff\left(\rt[[[][]]]\right)+ \coeff\left(\rt[[[][][][]]]\right)\right) , \\
            &1 + 12 \coeff\left(\rt[[[]][]]\right)+ 36 \coeff\left(\rt[[[][]][][]]\right)+ 12 \coeff\left(\rt[[[][][]]]\right)= 6 \left(6 \coeff\left(\rt[[[]][][]]\right)+ \coeff\left(\rt[[[][]]]\right)+ 4 \coeff\left(\rt[[[][][]][]]\right)\right) , \\
            &-1 + 6 \coeff\left(\rt[[[]]]\right)+ 6 \coeff\left(\rt[[[]][]]\right)- 3 \coeff\left(\rt[[[][]]]\right)+ 12 \coeff \left(\rt[[[[]]]]\right)- 36 \coeff\left(\rt[[[[]]][]]\right)- 18 \coeff\left(\rt[[[[][]]]]\right)+ 36 \coeff\left(\rt[[[[][]]][]]\right) = 0,\\
            &-(5/144) - \coeff\left(\rt[[[]]]\right)/6 + 2 \coeff\left(\rt[[[]][]]\right)/3 - \coeff\left(\rt[[[]][[]]]\right)/2 + \coeff\left(\rt[[[][]]]\right)/6 + \coeff\left(\rt[[[[]]]]\right)/6 - \coeff\left(\rt[[[[]][]]]\right)+ \coeff\left(\rt[[[[]][[]]]]\right) = 0, \\
            &\coeff \left(\rt[[][][][][]] \right) = \frac{1}{6}, \\
            &8 \left(\coeff\left(\rt[[[]]]\right)+ 3 (\coeff\left(\rt[[[]][[]]]\right)+ \coeff\left(\rt[[[[][]]]]\right)+ 2 \coeff\left(\rt[[[[[]]][]]]\right))\right) = 1 + 12 \coeff\left(\rt[[[][]]]\right)+ 8 \coeff\left(\rt[[[[]]]]\right)+ 48 \coeff\left(\rt[[[[]]][[]]]\right)+ 48 \coeff\left(\rt[[[[[]][]]]]\right), \\
            &2 \left(6 \coeff\left(\rt[[[]][]]\right)+ 12 \coeff\left(\rt[[[]][[]]]\right)+ 9 \coeff\left(\rt[[[][]]]\right)+ 8 \coeff\left(\rt[[[[]]]]\right)+ 24 \coeff\left(\rt[[[[]][]][]]\right)+     12 \coeff\left(\rt[[[[][]][]]]\right)\right) = 1 + 4 \coeff\left(\rt[[[]]]\right)+ 24 \coeff\left(\rt[[[][]][]]\right)+ 24 \coeff\left(\rt[[[][]][][]]\right)+ 24 \coeff\left(\rt[[[[]]][]]\right)\\
            &\qquad\phantom{2 \left(6 \coeff\left(\rt[[[]][]]\right)+ 12 \coeff\left(\rt[[[]][[]]]\right)+ 9 \coeff\left(\rt[[[][]]]\right)+ 8 \coeff\left(\rt[[[[]]]]\right)+ 24 \coeff\left(\rt[[[[]][]][]]\right)+     12 \coeff\left(\rt[[[[][]][]]]\right)\right) =} + 48 \coeff\left(\rt[[[[]][]]]\right)+   12 \coeff\left(\rt[[[[][]]]]\right), \\
            &-\coeff\left(\rt[[[]]]\right)^2 + \coeff\left(\rt[[[[]]]]\right)- 2 \coeff\left(\rt[[[[[]]]]]\right)+ 2 \coeff\left(\rt[[[[[[]]]]]]\right) = 0, \\
          \end{aligned}
          \endgroup
    \label{cond:p2q6}
    \end{equation}

Finally, using the standard expressions for the elementary differentials, we obtain the
following conditions on the Butcher coefficients.

For PEP order 3:

\begin{equation} \label{PEP3}
    \sum_i b_i c_i^2 = \frac{1}{3},
\end{equation}
For PEP order 4:
\begin{equation} \label{PEP4}
        \begin{aligned}
            \displaystyle \sum_{i,j,k}b_ia_{ij}a_{jk}c_k & = \displaystyle \sum_{i,j}b_ia_{ij}c_j -\frac{1}{8}, \\
            \displaystyle \sum_{i,j}b_ic_ia_{ij}c_j -\frac{1}{2} \displaystyle \sum_{i,j}b_ia_{ij}c_j^2 & = \frac{1}{12}, \quad &
            \sum_i b_i c_i^3 & = \frac{1}{4}. \\
        \end{aligned}
    \end{equation}
For PEP order 5:
\begin{equation} \label{PEP5}
        \begin{aligned}
            \displaystyle \sum_{i,j,k}b_ic_ia_{ij}a_{jk}c_k + \frac{1}{2}\displaystyle \sum_{i,j,k}b_ia_{ij}a_{jk}c_k^2 & = \displaystyle \sum_{i,j,k}b_ia_{ij}a_{jk}c_k &+ \frac{1}{2}\sum_{i,j}b_ia_{ij}c_j^2 - \frac{1}{2}\displaystyle \sum_{i,j}b_ia_{ij}c_j+\frac{1}{24}, \\
            2\displaystyle \sum_{i,j,k}b_ia_{ij}c_ja_{jk}c_k - \displaystyle \sum_{i,j}b_i(a_{ij}c_j)^2 & = \displaystyle \sum_{i,j,k}b_ia_{ij}a_{jk}c_k &+\displaystyle \sum_{i,j}b_ia_{ij}c_j^2 - \displaystyle \sum_{i,j}b_ia_{ij}c_j+\frac{1}{24} , \\
            \displaystyle \sum_{i,j}b_ic_i^2a_{ij}c_j - \frac{1}{3}\displaystyle \sum_{i,j}b_ia_{ij}c_j^3 & = \frac{1}{12}, \quad &
            \displaystyle \sum_ib_ic_i^4 & = \frac{1}{5}. \\
        \end{aligned}
    \end{equation}
To save space, we omit the final expressions for PEP order 6.

In Table~\ref{tab:orderconditions}, we compare the number of order
conditions for PEP and pseudo-symplectic methods.
\begin{table}[htbp]
\centering
\caption{Number of conditions for pseudo-symplectic and pseudo-energy-preserving methods.}
\label{tab:orderconditions}
\begin{tabular}{||ccc||}
\hline
Order        & Pseudo-Symplectic & Pseudo-Energy-Preserving  \\
        \hline
1 & 1       &    1      \\
\hline
2 & 2       & 2          \\
\hline
3  & 3       & 3        \\
\hline
4  & 6        & 6
\\ \hline
5  & 12       & 10       \\
\hline
6  & 28       & 21        \\
\hline
\end{tabular}
\end{table}

\begin{remark}
\label{rem:Hamiltonian_redundancies}
  Typically, in order to show that results on Runge-Kutta methods derived from a B-series analysis
  are necessary conditions (not just sufficient), one needs to use the fact that the elementary
  differentials are independent \cite[Section~314]{Butcher2003}.  The special
  structure of Hamiltonian systems \eqref{hamiltonian-system} raises the question
  whether the elementary differentials are still independent (enough) when the system is assumed
  to be Hamiltonian. This has been analyzed,
  e.g., in \cite{abia1993partitioned,calvo1994canonical},
  leading to sufficient and necessary conditions for symplecticity.
  In the case considered in this article, we note that the condition of
  Theorem~\ref{thm:energy-preserving-flow} of \cite{celledoni2010energy} is
  both necessary and sufficient for energy preservation of the flow.
\end{remark}

\section{Methods with maximal pseudo-energy-preserving order} \label{sec:methods}

Clearly, any pseudo-energy-preserving RK method must have $q \ge p$; i.e., the order of
energy preservation $q$ is at least the classical order $p$.  Since there is only one
condition for order two, and it is also the condition for energy preservation of
order two, we see that every method with $q=2$ must also have $p=2$.  Thus, the first
interesting class of methods to consider is that of $(s,p,q)=(2,2,3)$,
where $s \ge p$ is the number of stages of the explicit Runge-Kutta methods.
All two-stage, second-order explicit RK methods have the form
\begin{equation}
\renewcommand\arraystretch{1.2}\label{Met:rk22}
\begin{array}
{c|cc}
0 & 0 & 0\\
\alpha & \alpha & 0\\
\hline
& 1- \frac{1}{2\alpha} & \frac{1}{2\alpha}
\end{array}.
\end{equation}
The unique choice that yields $q=3$ is $\alpha = 2/3$ (see \cite{celledoni2009energy}):
\begin{align*}
\renewcommand\arraystretch{1.2}
\begin{array}
{c|cc}
0\\
\frac{2}{3} & \frac{2}{3}\\
\hline
& \frac{1}{4} &\frac{3}{4}.
\end{array}
\end{align*}

Similarly, among three-stage second-order explicit RK methods only
one satisfies the three conditions \eqref{PEP4} for PEP order four
(see \cite{celledoni2009energy}):
\begin{align*}
\renewcommand\arraystretch{1.2}
\begin{array}
{c|ccc}
0\\
\frac{1}{3} & \frac{1}{3}\\
\frac{5}{6}&  -\frac{5}{48} & \frac{15}{16}\\
\hline
& \frac{1}{10} &\frac{1}{2} &\frac{2}{5}.
\end{array}
\end{align*}

To find higher-order PEP methods, we have conducted a numerical search
using MATLAB's \verb|fmincon| constrained optimization routine (via the package RK-Opt \cite{ketcheson2020RK-Opt}),  with constraints
corresponding to the traditional order conditions as well as the PEP conditions \eqref{PEP3}-\eqref{PEP5}.
For each combination of stage count $s$ and classical order $p$, we have searched for
methods with the highest feasible value of the PEP order $q$.
In Table \ref{tab:maxepo} we list the highest PEP order $q$ of methods found
with a prescribed number of stages $s$ and classical order $p$.
We see that, except for the case $s=p=2$, additional stages are required
in order to achieve $q>p$.

Since \verb|fmincon| employs local search algorithms, and the problem is non-convex, we
cannot guarantee that the methods found are the best possible.  However, for each choice
of $(s,p,q)$ these are the results of many thousands of searches starting from different
randomly-generated initial guesses, enforcing the classical and PEP order
conditions, and minimizing the Euclidean norm of the leading truncation error coefficients.
We have chosen to show results only up to $s=5$ as the
optimization algorithms seem to
struggle greatly for larger stage numbers.  However, we have also found $q=6$ methods with
$(s,p)=(6,3)$ and $(s,p)=(7,4)$, which are included in some of the numerical examples later.

\begin{table}[htbp]
\centering
\caption{Maximal pseudo-energy-preserving order for given classical order $p$ and number of stages $s$ for explicit Runge-Kutta methods found in a numerical search.}
\label{tab:maxepo}
\begin{tabular}{||ccccc||}
\hline
        & $p = 2$ & $p = 3$ & $p = 4$ & $p = 5$ \\
        \hline
$s = 2$ & 3        & -       & -       & -       \\ \hline
$s = 3$ & 4        & 3       & -       & -      \\ \hline
$s=4$   & 5        & 4       & 4       & -       \\ \hline
$s=5$   & 6        & 5       & 5      & - \\ \hline
\end{tabular}
\end{table}

\section{Numerical tests} \label{sec:tests}
In this section we compare various Runge-Kutta methods, including:
\begin{itemize}
    \item Some classical methods, denoted by RK($s,p$); specifically, we use
          \begin{itemize}
            \item RK(2,2): the second-order explicit midpoint method of Runge \cite{runge1895numerische}
            \item RK(4,4): the classical fourth-order Runge-Kutta method of Kutta \cite{kutta1901beitrag}
          \end{itemize}
    \item Some pseudo-symplectic methods from \cite{aubry1998pseudo}, denoted by PS($s,p,q$);
    \item New pseudo-energy-preserving methods found herein by numerical search, denoted by PEP($s,p,q$).
\end{itemize}
Here $s, p$ denote the number of stages and classical order, while $q$ denotes the
order of pseudo-symplecticity or energy preservation.  We focus on application to Hamiltonian systems
and we are primarily interested in understanding the long-time accuracy of these methods.
It is known that exactly energy-preserving numerical methods exhibit linear error growth (versus
quadratic error growth for generic methods) in many cases; see e.g. \cite{duran1998numerical,duran2000numerical,ranocha2021rate}.
It is natural to expect that PEP methods
(and possibly PS methods) will exhibit linear error growth for moderately long times,
since the dominant error terms do not involve a change in the energy.


The theory of Theorem \ref{thm:energy-preserving-flow} is established exclusively for canonical Hamiltonian systems, so we have dedicated Section \ref{sec:canonical} to
compare PEP methods with classical Runge-Kutta methods and PS methods for different canonical Hamiltonian systems, while in  Section \ref{sec:noncanonical} we test how PEP methods
behave for non-canonical Hamiltonian systems.

Since the PEP and PS methods often have a larger number of stages, and thus a larger cost per step,
we use differing step sizes in order to compare the methods fairly.
For each test, we use stepsizes given by
\begin{align}
    h = s \stepsizex,
\end{align}
where $s$ is the number of stages and the value $\stepsizex$ (the time
advancement per stage) is the same for all methods. The value of $\stepsizex$
for each test is given in the description of the corresponding figure.

 \subsection{Canonical Hamiltonian systems}
 \label{sec:canonical}

\subsubsection{Exponential entropy system}
The exponential entropy system
\begin{equation}
    \frac{d}{dt} \binom{u_1(t)}{u_2(t)} = \binom{-\exp(u_2(t))}{\exp(u_1(t))}, \quad u^0 = \binom{1}{0.5},
    \label{eq:exp_entropy_system}
\end{equation}
with Hamiltonian
\begin{equation}
    H(u) = \exp(u_1)+\exp(u_2), \quad \nabla_u H = \binom{\exp(u_1)}{\exp (u_2)},
\end{equation}
was studied for relaxation Runge-Kutta methods in \cite{Ranocha2020Entropy}.

We use this system to check the theoretical convergence rates of the proposed methods.
For this purpose, we start by solving  \eqref{eq:exp_entropy_system}
and calculate the error at $T=160$ for both the solution and the energy.
Let $e_i, e_{i+1}$ denote the errors obtained using stepsizes $h_i,h_{i+1}$ respectively.
Then we define the experimental order of convergence (EOC)
\begin{equation}
    EOC = \frac{\log( e_{i+1} / e_{i} )}{\log( h_{i+1} / h_{i} )} = \frac{\log( e_{i+1} / e_{i} )}{\log(2)}.
\end{equation}

\begin{table}[htbp]
\centering
  \caption{Convergence test for the exponential entropy system using PEP(6,3,6).}
\label{tab:EOC_PEP636}
\begin{tabular}{|l|c|c|c|c|}
\hline
$h$ & Solution Error & Solution EOC & Energy Error & Energy EOC \\
\hline
1/2      & 1.93e-01 & 5.05 & 1.06e-03 & 5.96 \\
1/4     & 5.81e-03 & 3.68 & 1.70e-05 & 5.62 \\
1/8    & 4.53e-04 & 3.14 & 3.47e-07 & 5.83 \\
1/16   & 5.15e-05 & 3.01 & 6.08e-09 & 5.92 \\
1/32  & 6.39e-06 & 3.00 & 1.00e-10 & 5.96 \\
1/64 & 8.00e-07 & 2.99 & 1.61e-12 & 5.36 \\
\hline
\end{tabular}
\end{table}

\begin{table}[htbp]
\centering
  \caption{Convergence test for the exponential entropy system using PEP(7,4,6).}
\label{tab:EOC_PEP746}
\begin{tabular}{|c|c|c|c|c|}
\hline
$h$ & Solution Error & Solution EOC & Energy Error & Energy EOC \\
\hline
1/2      & 5.84e-01 & 6.59 & 3.62e-03 & 6.68 \\
1/4      & 6.05e-03 & 6.56 & 3.54e-05 & 7.25 \\
1/8      & 6.40e-05 & 5.02 & 2.32e-07 & 7.79 \\
1/16     & 1.97e-06 & 4.08 & 1.05e-09 & 8.13 \\
1/32     & 1.16e-07 & 3.96 & 3.74e-12 & 4.19 \\
1/64     & 7.50e-09 & 5.95 & 2.05e-13 & 2.64 \\
\hline
\end{tabular}
\end{table}

Errors and convergence rates are displayed in Tables \ref{tab:EOC_PEP636} and \ref{tab:EOC_PEP746}.
For the PEP(6,3,6) method we see convergence rates that closely match the
theoretical values $p=3, q=6$, until rounding errors begin to affect the results at smaller step sizes.
For the PEP(7,4,6) method the results are less clear, but roughly consistent with the theory,
as we see a consistently higher rate of convergence for the energy compared to the solution.
Convergence tests for the other canonical systems studied below are not shown but gave results
similar to those displayed here.

\subsubsection{Undamped Duffing Oscillator}

The equation for the undamped Duffing oscillator is
\begin{equation}
    \frac{d^2}{dt^2}u_1(t) = u_1(t) - u_1(t)^3,
    \label{eq:duffing_oscillator}
\end{equation}
and its energy is given by
\begin{equation}
    H(u_1,u_2) = \frac{1}{2}u_2^2 -\frac{1}{2}u_1^2+\frac{1}{4} u_1^{4},
\end{equation}
where $u_2$ is the momentum given by $u_2(t) = \frac{d}{dt}u_1(t)$.
We take initial data $(u_1,u_2) = (1.4142, 0)$.

\begin{figure}[htbp]
\centering
\begin{subfigure}[b]{0.49\textwidth}
  \centering
  \includegraphics[width=\textwidth]{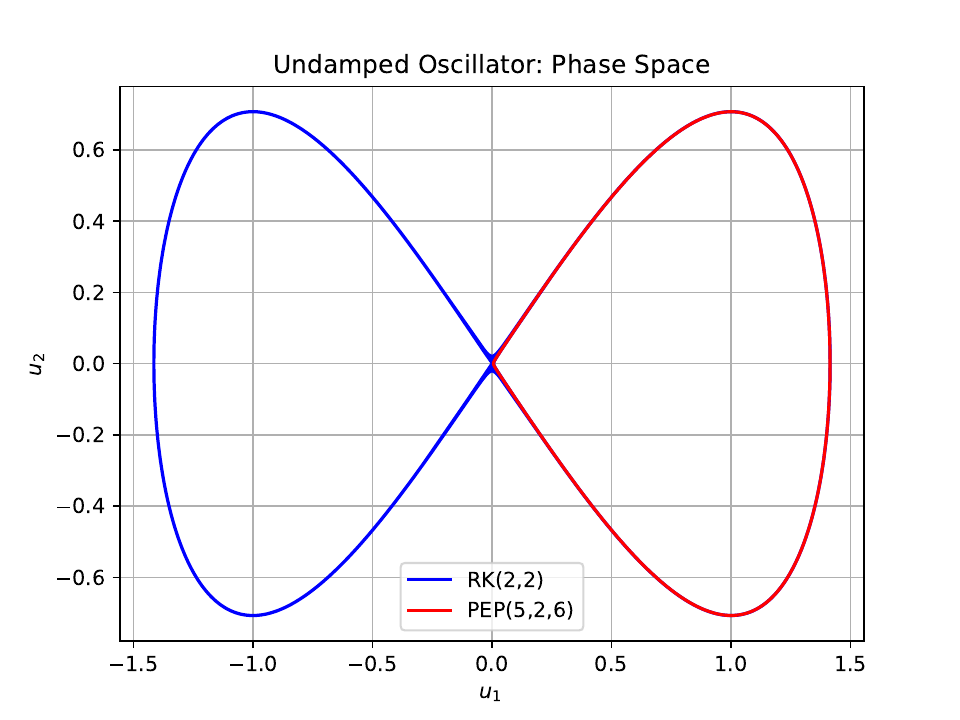}
  \caption{RK(2,2)  in contrast to PEP(5,2,6) with $\stepsizex = 1/200$.}
  \label{fig:undampedrk2}
\end{subfigure}%
\hspace*{\fill}
\begin{subfigure}[b]{0.49\textwidth}
  \centering
  \includegraphics[width=\textwidth]{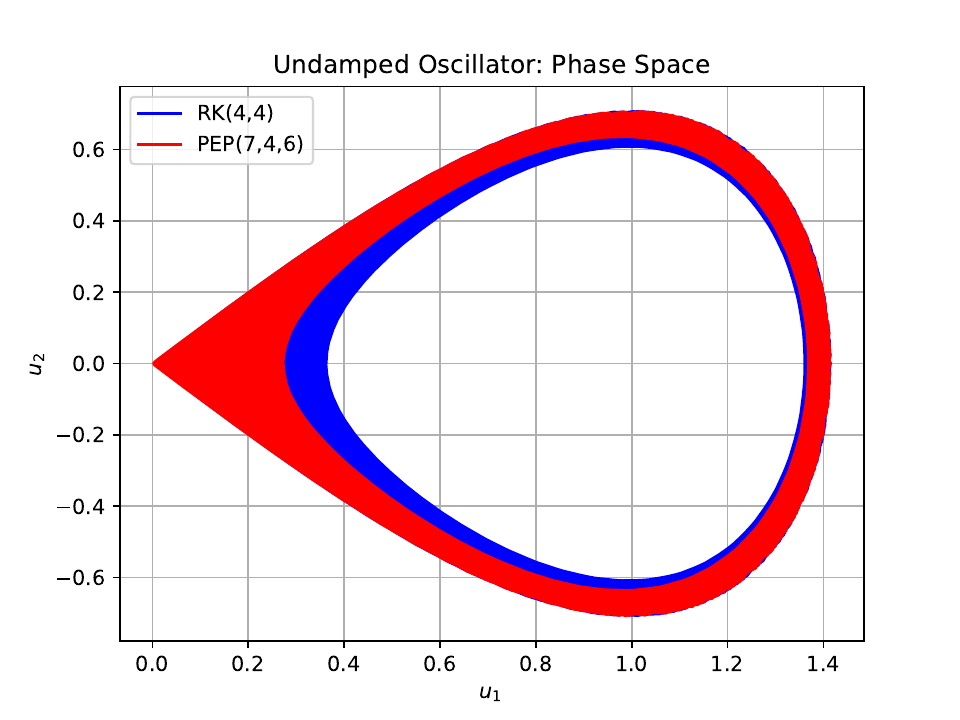}
  \caption{RK(4,4)  in contrast to PEP(7,4,6) with $\stepsizex = 1/20$.}
  \label{fig:undampedrk4}
\end{subfigure}
\caption{Simulation of the undamped duffing oscillator comparing RK methods for which $p = q$ with RK methods with higher PEP.}
\label{fig:undampedpep}
\end{figure}

This initial value lies on a periodic orbit just inside a homoclinic connection, and the correct solution is shown in red in Figure 6(a). Numerical errors can cause two types of behavior, depending on the sign of the energy error: If the energy increases, the solution can approximate the mirror-image homoclinic orbit, which lies in the left half-plane, as shown in Figure \ref{fig:undampedrk2}. If the energy decreases, the solution can move inward toward the center equilibrium, as shown in Figure \ref{fig:undampedrk4}.

For methods of order 2 and 3, we typically see the energy increase, but PEP methods are able to remain near the correct (right half-plane) orbit under much larger step sizes. For example, letting $h_{max}$
denote the largest step size for which the orbit stays in the right half-plane, we find $h_{max} = 0.004$ for RK(2,2) and $h_{max} = 0.152$ for PEP(5,2,6); see Figure \ref{fig:undampedrk2}.
For the 4th-order methods we tested, the energy decreases over time but (as expected) does so more slowly for PEP methods, as illustrated in Figure \ref{fig:undampedrk4}.

\subsubsection{Hénon-Heiles System}
\label{sec:henon-heiles}

\begin{figure}[htbp]
  \centering
  \begin{subfigure}[b]{0.49\textwidth}
    \centering
    \includegraphics[width=\textwidth]{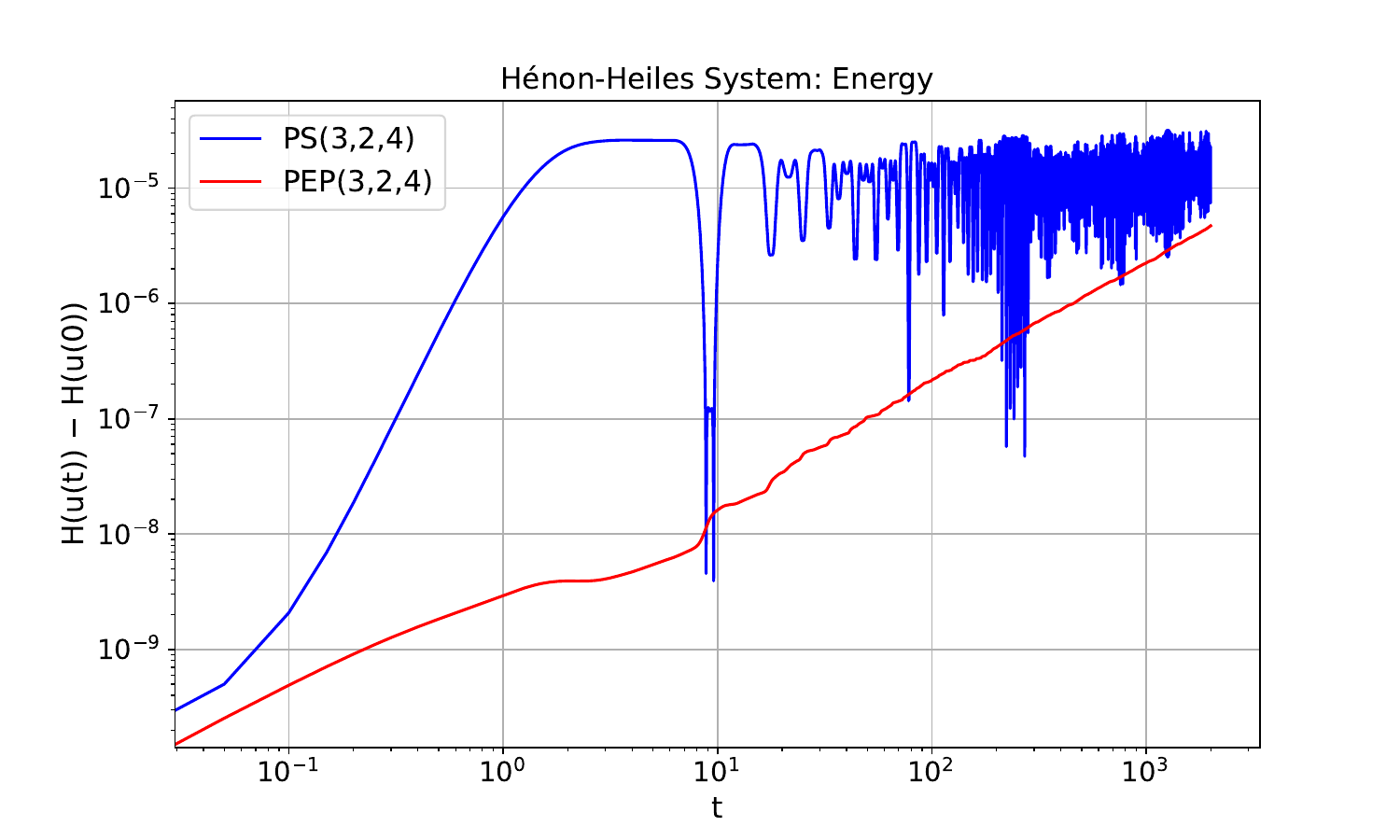}
    \caption{PS(3,2,4)  in contrast to PEP(3,2,4) with $\stepsizex = 5/300$.}
    \label{fig:hh324}
  \end{subfigure}%
  \hspace*{\fill}
  \begin{subfigure}[b]{0.49\textwidth}
    \centering
    \includegraphics[width=\textwidth]{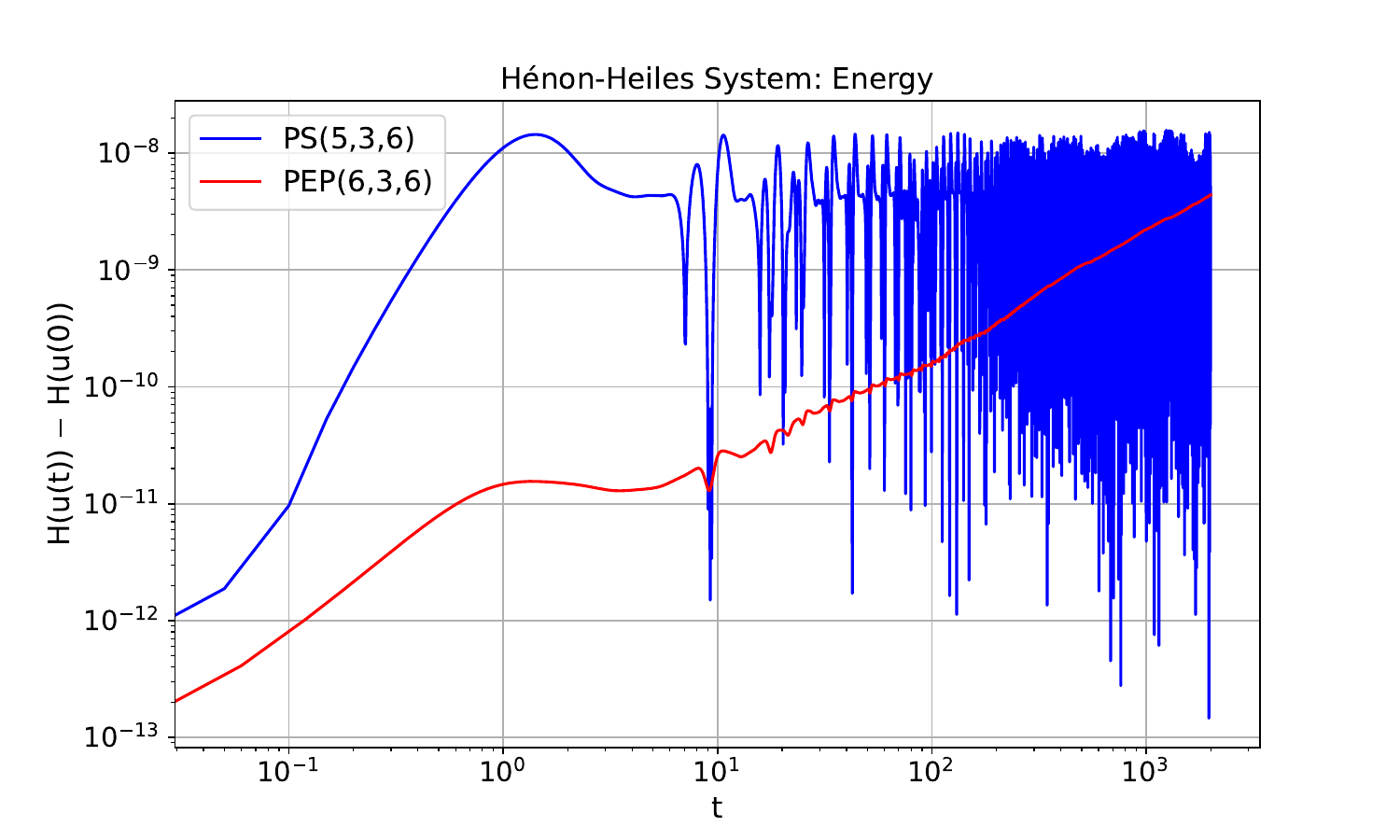}
    \caption{PS(5,3,6)  in contrast to PEP(6,3,6)  with $\stepsizex = 1/100$.}
    \label{fig:hh636}
  \end{subfigure}
  \caption{Numerical energy evolution for the Hénon-Heiles system, comparing PS and PEP methods.}
  \label{fig:hhsystem}
  \end{figure}

Another ilustrative comparison between classical Runge-Kutta methods and energy-preserving methods is given by analyzing the Hénon-Heiles system, described by the Hamiltonian
\begin{equation}
  H = \frac{1}{2}\left( p_x^2 + p_y^2 \right) + \left( x^2 + y^2 \right) + \lambda \left( x^2y- \frac{y^3}{3} \right).
  \label{eq:hhsystem}
\end{equation}
Following \cite{quispel2008newclass}, we consider the initial condition
$$
x = 0.1, \; \; y = -0.5, \; \; \dot{x} = 0, \; \; \dot{y} = 0
$$
and we set $\lambda=1$.
In this case the exact solution remains within a certain triangular region,
as discussed and tested in \cite{quispel2008newclass}.
For all of the methods considered here, the numerical solution remains
in this invariant region for rather long times, so it is less useful as a metric for comparison.
Instead, we focus on the value of $H$ over time obtained with PEP and PS methods,
as shown in Figure \ref{fig:hhsystem}.  We see that PEP methods provide significantly
better conservation of energy.

\subsubsection{Benjamin-Bona-Mahony Equation}

As a final example, we consider the Benjamin-Bona-Mahony (BBM) equation
\cite{benjamin1972model}
\begin{equation}
    u_t + u_x +uu_x - u_{txx} = 0
    \label{eq:bbm}
\end{equation}
with periodic boundary conditions on the domain
$[-90, 90]$. We take a solitary wave as the initial condition:
\begin{equation}
  u(t, x) = \frac{A}{\cosh\bigl(K (x - c t) \bigr)^2},
  \quad
  \text{with }
  A = 3 (c - 1),
  \;
  K = \frac{1}{2} \sqrt{1 - 1 / c},
  \;
  c = 1.2.
\end{equation}
We discretize in space using Fourier pseudospectral collocation, using the Julia
\cite{bezanson2017julia} package SummationByPartsOperators.jl
\cite{ranocha2021sbp}, wrapping FFTW.jl \cite{frigo2005design}.
We use the following conservative semidiscretization
given in \cite{linders2023iterativemethods}:
\begin{equation} \label{semidisc}
    u_t + (I -D_2)^{-1}D_1\left( \frac{1}{2}u^2 +u\right) = 0,
\end{equation}
Here $D_1$ and $D_2$ are the Fourier spectral collocation
differentiation matrices of order 1 and 2, respectively.
The semi-discrete scheme conserves a discrete mass and
a discrete approximation of the Hamiltonian
\begin{equation}
  \int \left( \frac{1}{2} u^2 + \frac{1}{6} u^3 \right),
\end{equation}
as discussed in \cite{linders2023iterativemethods}.

\begin{figure}
    \centering
    \includegraphics[scale = 0.6]{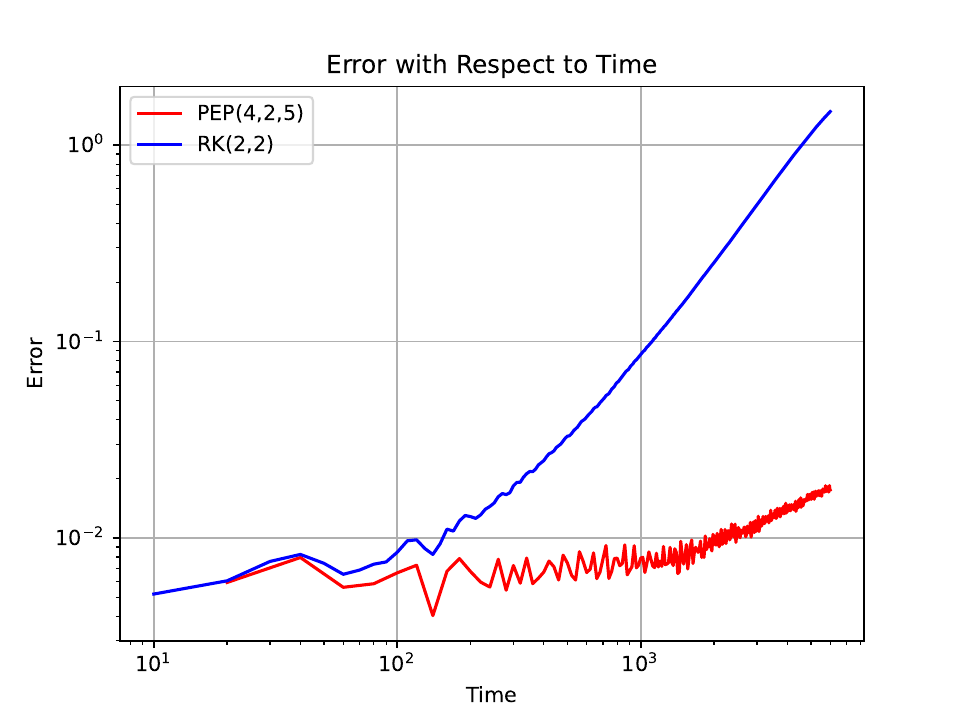}
    \caption{Error growth for the BBM equation comparing the midpoint method with with  PEP(4,2,5) with $\stepsizex = 0.05$.}
    \label{fig:BBM}
\end{figure}

We study the behavior of solutions when \eqref{semidisc} is
solved with two different time integration methods: a traditional
RK(2,2) method with $\alpha=\frac{1}{2}$, and our PEP(4,2,5) method.
The error over time is shown in Figure \ref{fig:BBM}; we see
again that the method with higher PEP order produces smaller errors
over long times.

\subsection{Non-canonical Hamiltonian systems}
\label{sec:noncanonical}

Theorem \ref{thm:energy-preserving-flow} applies only to canonical Hamiltonian systems, for which the skew-symmetric
matrix $J$ is constant.  Here we test whether any benefit is gained by using PEP methods for
non-canonical Hamiltonian systems, in which $J$ depends on the solution $u$.

\subsubsection{Nonlinear oscillator}

Consider the simple nonlinear oscillator problem
\begin{equation}
    \frac{d}{dt}\binom{u_1(t)}{u_2(t)} = ||u(t)||^{-2}\binom{-u_2(t)}{u_1(t)}, \quad u^0 = \binom{1}{0},
    \label{eq:nonlinearoscillator}
\end{equation}
with conserved energy $\|u\|_2^2$, which has been studied for example in \cite{ranocha2021strong,Ranocha2020RRKHamiltonian}.

\begin{figure}[htbp]
\centering
\begin{subfigure}[b]{0.49\textwidth}
  \centering
  \includegraphics[width=\textwidth]{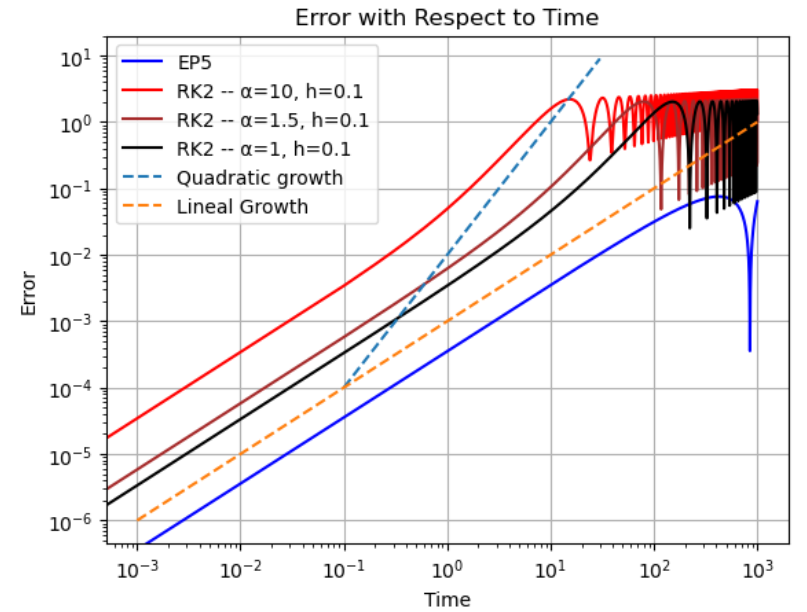}
  \caption{Comparison of 2nd-order methods on problem \eqref{eq:nonlinearoscillator}, with $\stepsizex = 1/20$.}
  \label{fig:nonlinearoscillatorrk2}
\end{subfigure}%
\hspace*{\fill}
\begin{subfigure}[b]{0.49\textwidth}
  \centering
  \includegraphics[width=\textwidth]{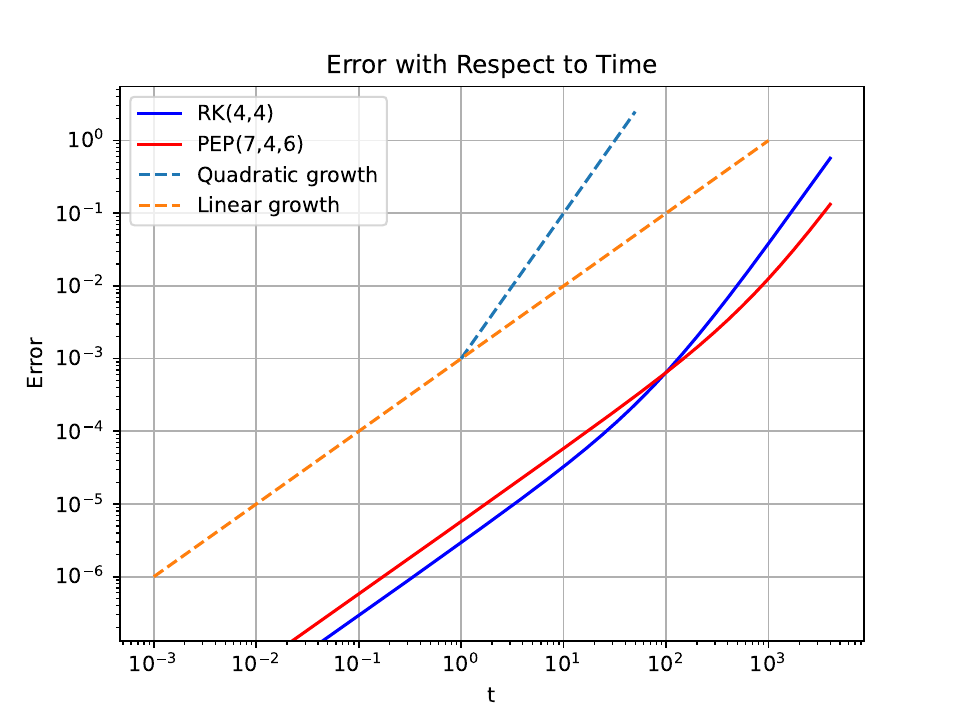}
  \caption{ RK(4,4)  in contrast to PEP(7,4,6), with $\stepsizex = 1/40$.}
  \label{fig:nonlinearoscillatorrk4}
\end{subfigure}
\caption{Error growth for the nonlinear oscillator problem
\eqref{eq:nonlinearoscillator}, comparing RK methods for which $p = q$ with RK
methods with higher PEP order. } \label{fig:nonlinearoscillator}
\end{figure}

We solve this problem with a range of methods with different stepsizes each and plot the error as a function of
time in Figure~\ref{fig:nonlinearoscillatorrk2}.
We include some classical two-stage second-order methods of the form \eqref{Met:rk22}
with different values of $\alpha$, as well as a PEP(4,2,5) method.
Note that for $\alpha=2/3$ we recover the PEP(2,2,3) method.  The error growth is initially linear, but eventually
becomes superlinear before saturating (due to the periodicity of the solution).
We see that as $\alpha$ gets closer to the optimal $2/3$
value, the error growth remains linear up to later times.
The PEP(4,2,5) method exhibits linear error growth over a longer time than
any of the methods with $q=2$.

In Figure \ref{fig:nonlinearoscillatorrk4} we compare the classical 4th-order
method to a 4th-order PEP method with $q=6$.  We see that the latter method
gives a larger error initially but smaller error for long times, since the
error growth remains linear until later.

\begin{figure}[htbp]
\centering
\begin{subfigure}[b]{0.49\textwidth}
  \centering
  \includegraphics[width=\textwidth]{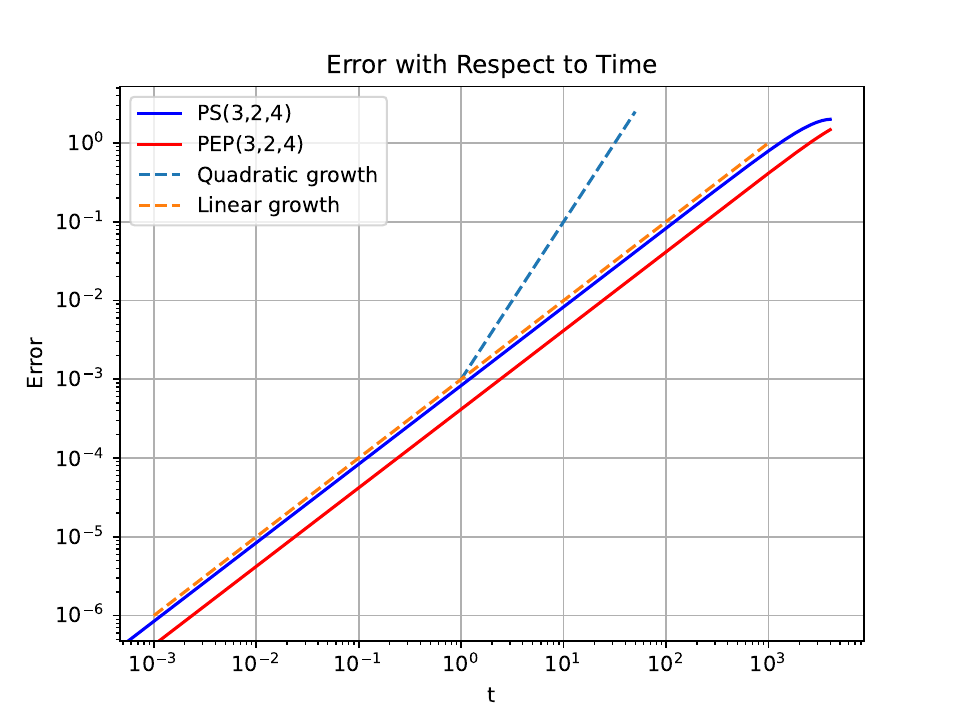}
  \caption{PS(3,2,4) vs. PEP(3,2,4), with $\stepsizex = 1/30$.}
  \label{fig:nonlinearsym24}
\end{subfigure}%
\hspace{\fill}
\begin{subfigure}[b]{0.49\textwidth}
  \centering
  \includegraphics[width=\textwidth]{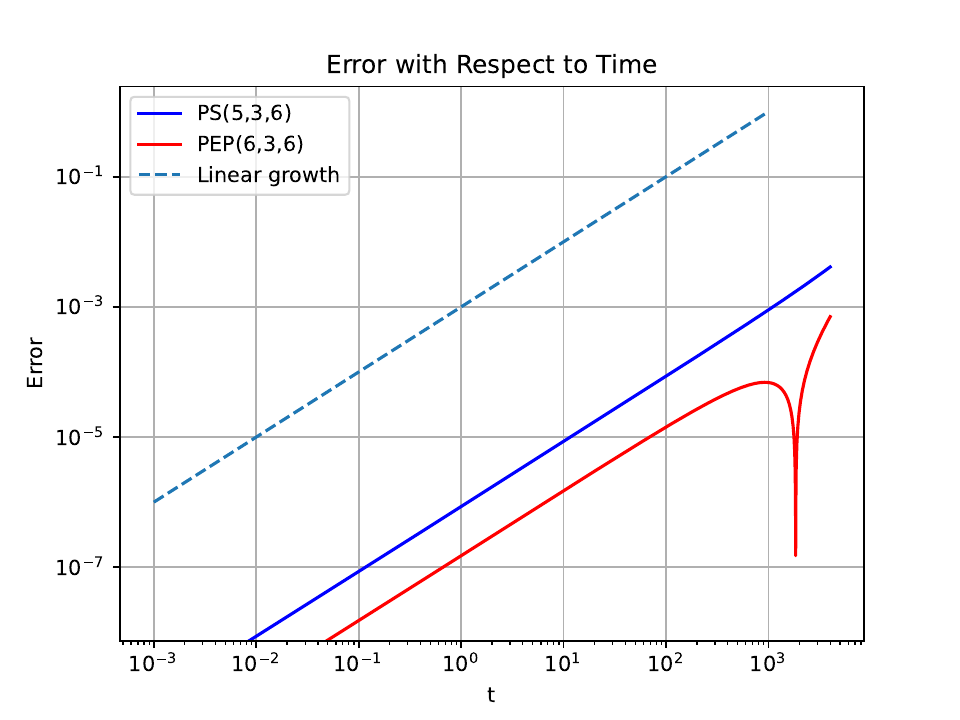}
  \caption{PS(5,3,6)  vs. PEP(6,3,6) with $\stepsizex = 1/50$. }
  \label{fig:nonlinearsym36}
\end{subfigure}
\caption{Error growth for the nonlinear oscillator comparing PS with PEP methods.}
\label{fig:nonlinear_sym}
\end{figure}

In Figure \ref{fig:nonlinear_sym} we compare PEP and pseudo-symplectic (PS) methods with the same orders and pseudo-orders.
We see that these methods behave very similarly, exhibiting linear error growth over the whole time range.

\subsubsection{Lotka-Volterra Equations}

The Lotka-Volterra system \cite{hairer2006geometric}
\begin{equation}
    \frac{d}{dt}\binom{u_1(t)}{u_2(t)} = \binom{u_1(t)(1- u_2(t))}{u_2(t)(u_1(t) -1 )}, \; \; u(0) = \binom{1}{2},
    \label{eq:lotkavolterra}
\end{equation}
is a model for the populations of two competing biological species.
The Hamiltonian $H(u)= u_1 + u_2 - \log(u_1) - \log(u_2)  $ is constant in time.

Results of a convergence test are displayed in Table \ref{tab:eoc_lv}.  We see that the energy EOC is close to the classical
order of the method, not its PEP order, reflecting the fact that PEP methods are based on analysis of canonical Hamiltonian
systems only.

\begin{figure}[htbp]
\centering
\begin{subfigure}[b]{0.49\textwidth}
  \centering
  \includegraphics[width=\textwidth]{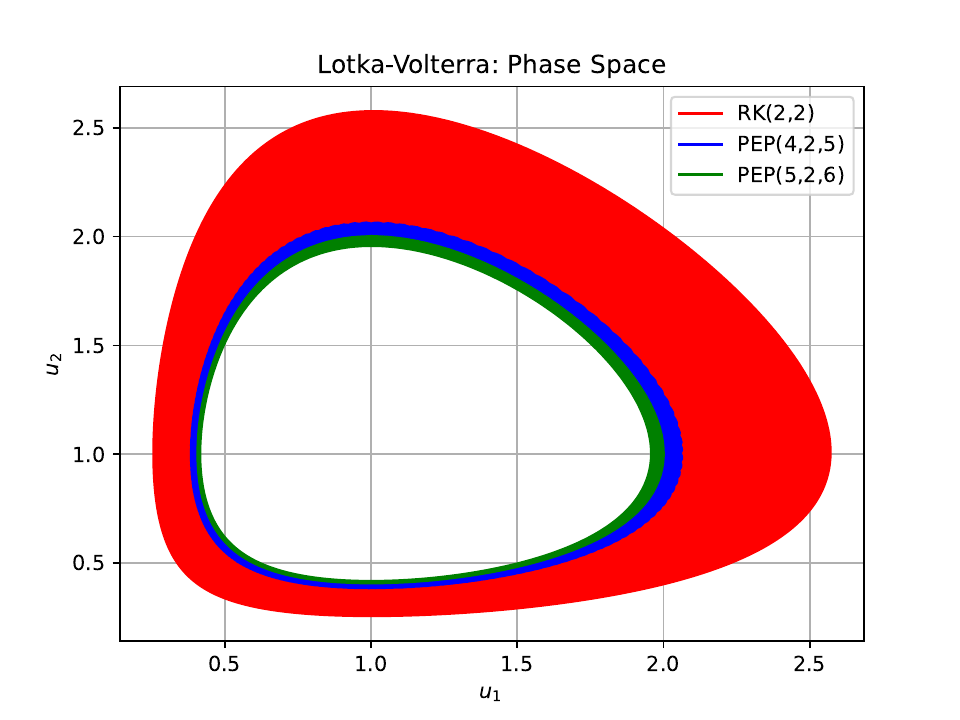}
  \caption{Lotka-Volterra phase space.}
  \label{fig:LVtest}
\end{subfigure}%
\hspace{\fill}
\begin{subfigure}[b]{0.49\textwidth}
  \centering
  \includegraphics[width=\textwidth]{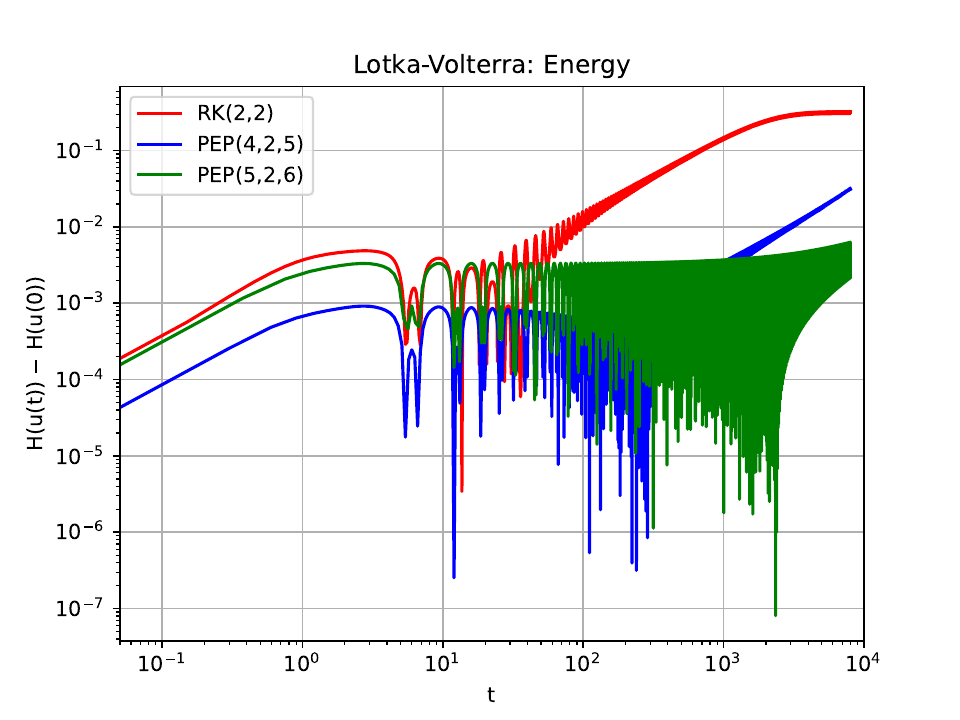}
  \caption{Lotka-Volterra energy.}
  \label{fig:LVenergy}
\end{subfigure}
\caption{Numerical solution for the Lotka-Volterra equations comparing the behavior of RK(2,2) with $\alpha = 0.5$, PEP(4,2,5)  and PEP(5,2,6), with $\stepsizex = \frac{3}{40}$.}
\label{fig:LV}
\end{figure}

\begin{figure}[htbp]
\centering
\begin{subfigure}[b]{0.49\textwidth}
  \centering
  \includegraphics[width=\textwidth]{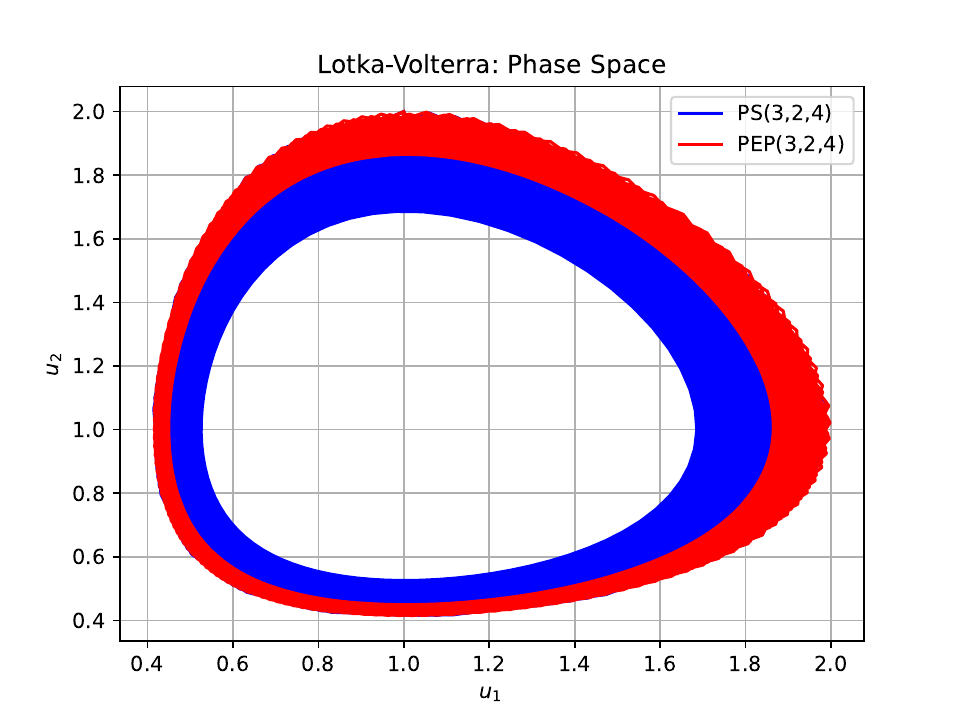}
  \caption{Lotka-Volterra phase space.}
  \label{fig:LVtestsym24}
\end{subfigure}%
\hspace{\fill}
\begin{subfigure}[b]{0.49\textwidth}
  \centering
  \includegraphics[width=\textwidth]{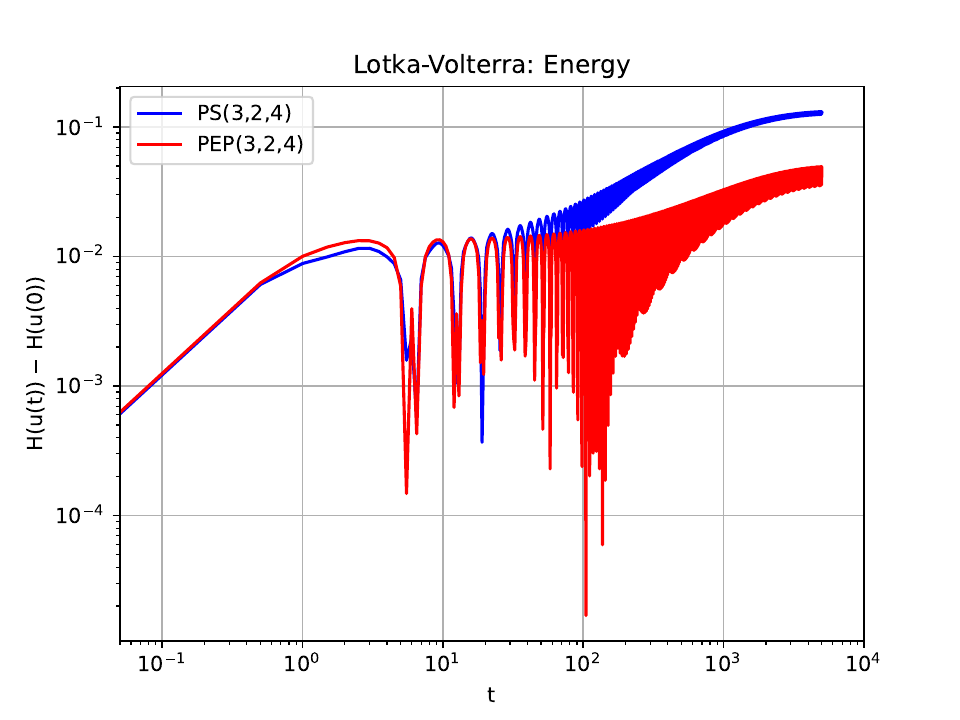}
  \caption{Lotka-Volterra energy.}
  \label{fig:LVenergysym24}
\end{subfigure}
\caption{ PS(3,2,4) vs. PEP(3,2,4), with $\stepsizex = 1/6$. }
\label{fig:LVsym24}
\end{figure}

\begin{figure}[htbp]
\centering
\begin{subfigure}[b]{0.49\textwidth}
  \centering
  \includegraphics[width=\textwidth]{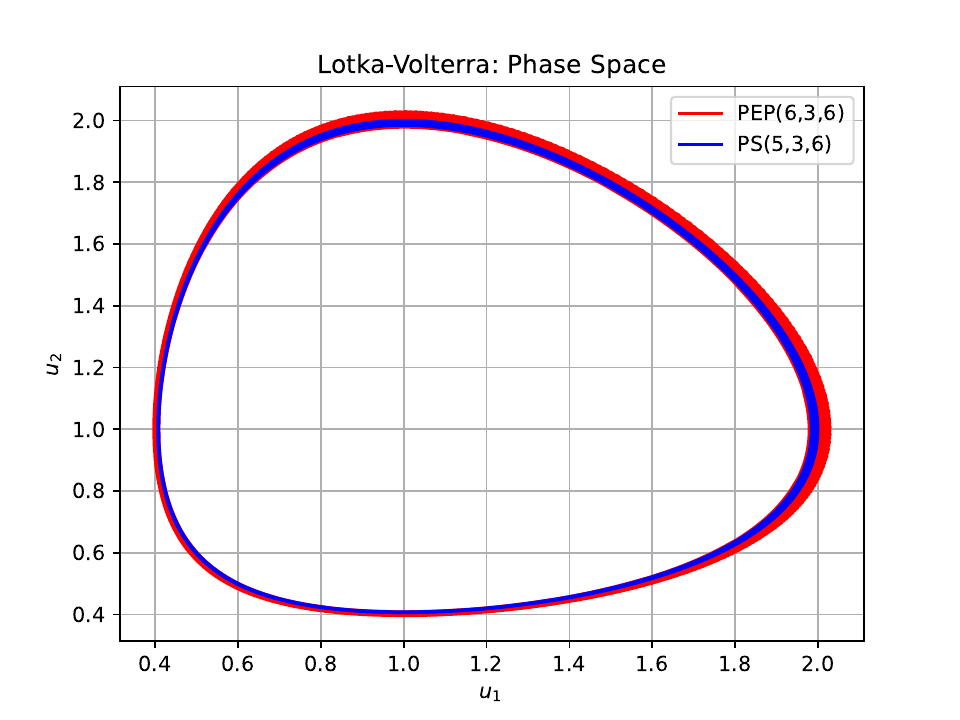}
  \caption{Lotka-Volterra phase space.}
  \label{fig:LVtestsym36}
\end{subfigure}%
\hspace{\fill}
\begin{subfigure}[b]{0.49\textwidth}
  \centering
  \includegraphics[width=\textwidth]{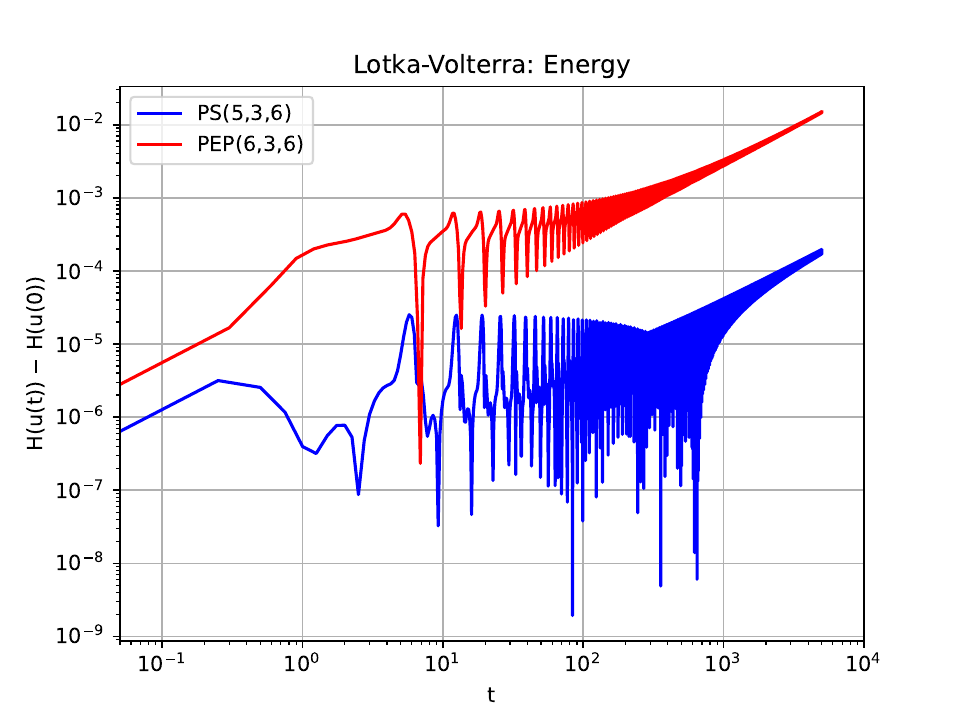}
  \caption{Lotka-Volterra energy.}
  \label{fig:LVenergysym36}
\end{subfigure}
\caption{ PS(5,3,6) vs. PEP(6,3,6) with $\stepsizex = 1/20$. }
\label{fig:LVsym36}
\end{figure}

Nevertheless, it seems that there is some benefit obtained by using a PEP method for this
non-canonical Hamiltonian system.
We compare the solutions and energy errors obtained with a variety of methods
in Figures \ref{fig:LV}, \ref{fig:LVsym24} and \ref{fig:LVsym36}.  We see that
methods with higher PEP order give smaller errors after long times, and stay visibly closer
to the correct periodic orbit.   Pseudo-symplectic methods show a similar benefit.

    \begin{table}[h!]
      \centering
      \caption{Convergence test for the Lotka-Volterra equations using PEP(6,3,6)}
      \label{tab:eoc_lv}
      \begin{tabular}{|c|c|c|}
      \hline
      $h$ & Energy Error & Energy EOC \\
      \hline
      1/2      & 7.73e-02 & 4.96e+00 \\
      1/4      & 2.48e-03 & 4.59e+00 \\
      1/8      & 1.03e-04 & 4.13e+00 \\
      1/16     & 5.88e-06 & 3.52e+00 \\
      1/32     & 5.12e-07 & 3.17e+00 \\
      1/64     & 5.70e-08 & 3.05e+00 \\
      \hline
      \end{tabular}
      \end{table}

For small values of $h$, both the PS and PEP methods in Figure \ref{fig:LVsym24} have a similar behavior. Then, we used a stepsize $h = 0.5$ in order to visualize the difference between both methods in Figure \ref{fig:LVsym24}.

\section*{Acknowledgements}

The first two authors were funded by the
King Abdullah University of Science and Technology (KAUST), and
through the KAUST Visiting Student Research Program.
The last author was supported by the Deutsche Forschungsgemeinschaft
(DFG, German Research Foundation, project number 513301895)
and the Daimler und Benz Stiftung (Daimler and Benz foundation,
project number 32-10/22).

\section*{Data availability statement}

The code used in this paper is available online in a Github repository:
\url{https://github.com/Sondar74/PEP_Reproducibility2024} \cite{barrios2024pseudoRepro}.

\appendix
\section{Method coefficients}
Here we list the coefficients of the PEP RK methods found and tested in this work.
We use the naming convention PEP($s$,$p$,$q$) where $s$ is the number of
stages, $p$ is the order of the method, and $q$ is the PEP order.

\bigskip

PEP(4,2,5):

\begin{equation}
\renewcommand\arraystretch{1.2}
\label{Met:EP5}
\begin{array}
{c|cccc}
0\\
\frac{1}{10} & \frac{1}{10}\\
\frac{19}{20} & -\frac{35816}{35721} &\frac{56795}{35721} \\
\frac{37}{63} & \frac{11994761}{5328000} & -\frac{11002961}{4420800} &  \frac{215846127}{181744000} \\
\hline
& -\frac{17}{222}  & \frac{6250}{15657} &\frac{5250987}{10382126} &\frac{4000}{23307}
\end{array}
\end{equation}

\begin{table}[h]
    \centering
    \caption{Floating point coefficients of the method PEP(5,2,6).}
    \label{Met_EP6}
    \begin{tabular}{|l|l|l|}
    \hline
    $a_{21} = 0.193445628056365 $ & $a_{31} = -0.090431947690469$ & $a_{32} = 0.646659568003039$ \\
    $a_{41} =-0.059239621354435$ & $a_{42} = 0.598571867726670$ & $a_{43} =  -0.010476084304794$ \\
    $a_{51} = 0.173154586278662$ & $a_{52} = 0.043637751980064$ & $a_{53} = 0.949323298732961$ \\
    $a_{54} =-0.262838451019868$ & $b_1 = 0.054828314201395$ & $b_2 = 0.310080077556546$\\
    $b_3 = 0.531276882919990$ & $b_4 = -0.135494569336049$ & $b_5 = 0.239309294658118$\\
    $c_2 = 0.193445628056365$ &$c_3 = 0.55622762031257$ & $c_4 = 0.528856162067441$\\
    $c_5 = 0.9032771859718189$ & & \\
    \hline
    \end{tabular}
\end{table}

A third-order method requires a bigger number of stages, according to Table \ref{tab:maxepo}. The following arrays are examples of RK methods of orders 3, 4 and 5:


\begin{table}[h]
    \centering
    \caption{Floating point coefficients of the method PEP(6,3,6).}
    \label{Met:pep36}
    \begin{tabular}{|l|l|l|}
    \hline
    $a_{21} = 0.12316523079127038$&  $a_{31} = -0.53348119048187126$&  $a_{32} = 1.1200645707708279$ \\
    $a_{41} = 0.35987162974687092$ &  $a_{42} = -0.17675778446586507$ &  $ a_{43} = 0.7331973326225617$\\
$    a_{51} = 0.015700424346522388$&$a_{52} = 0.02862938097533644$&  $ a_{53} = -0.014047147149911631$\\
$    a_{54} = -0.015653338246176568$ &   $a_{61} = -1.9608805853984794$ & $ a_{62} = -0.82154709029385564$\\
    $a_{63} = -0.0033631561953843502$ &  $a_{64} = 0.046367461001250457$& $ a_{65} = 2.782035718578454$\\
    $b_1 = 0.78642719559722885$&  $b_2 = 0.69510370728230297$& $b_3 = 0.42190724518033551$\\
    $b_4 = 0.21262030193155254$&  $b_5 = -0.70167978222250704$& $ b_6 = -0.41437866776891263$ \\
     $c_2 = 0.12316523079127044$  & $ c_3 = 0.58658338028895673$ &
    $c_4 = 0.91631117790356775$ \\
    $c_5 = 0.014629319925770667$  &  $c_6 = 0.042612347691984923$ & \\
    \hline
    \end{tabular}
\end{table}


\begin{table}[h]
    \centering
    \caption{Floating point coefficients of the method PEP(7,4,6).}
    \begin{tabular}{|l|l|l|}
    \hline
    $a_{21} = -0.10731260966924323$& $ a_{31} = 0.14772934954602848$& $ a_{32} = -0.12537555684690285$ \\
    $a_{41} = 0.7016079790308741$& $  a_{42} = -0.75094597518803941$& $ a_{43} = 0.76631666070124027$\\
    $a_{51} = -0.8967481787471202$& $ a_{52} = -0.43795858531068965$& $ a_{53} = 1.7727346351832869$\\
    $a_{54} = 0.1706052810617312$& $ a_{61} = 1.6243872270239892$& $ a_{62} = -0.69700589895015241$\\
    $a_{63} = -0.3861309831750398$& $ a_{64} = -0.032848941899304235$& $ a_{65} = 0.30227620385295728$\\
    $a_{71} = -0.32463926305048885$& $ a_{72} = -0.3480143346241919$& $ a_{73} = 1.3500419757109139$\\
    $a_{74} = 0.039096802121597336$& $ a_{75} = -0.17851883247877129$& $ a_{76} = 0.010142489530892661$\\
    $b_1 = -0.69203318482299292$& $ b_2 = 0.0074442860308153933$& $ b_3 = 0.93216717844052677$\\
    $b_4 = -1.159431111205361$&  $b_5 = 0.27787978605406632$&  $b_6 = 0.93890392164164138$\\
    $b_7 =0.69506912386130404$ &   $c_2 = -0.10731260966924323$& $ c_3 = 0.022353792699125609$\\
    $c_4 = 0.71697866454407488$& $ c_5 = 0.60863315218720804$& $ c_6 = 0.81067760685245005$\\
    $c_7 = 0.54810883720995185$ & & \\
    \hline
    \end{tabular}
\end{table}


\begin{table}[h]
    \centering
    \caption{Floating point coefficients of the method PEP(7,5,6).}
    \begin{tabular}{|l|l|l|}
    \hline
    $a_{21} = 0.34288981581855521$& $ a_{31} = 0.16800230418143236$& $ a_{32} = 0.1262987524809161$\\
    $a_{41} = 0.4326925567104672$& $  a_{42} = -0.24221982610439177$& $ a_{43} = 0.15241708521248304$\\$
    a_{51} = 0.019843989305203335$& $ a_{52} = 0.20330206481276515$& $ a_{53} = -0.3494376489494413$\\$
    a_{54} = 0.09780248603799992$& $ a_{61} = 3.5441758455721732$& $ a_{62} = 9.884560134482289$\\ $
    a_{63} = -3.7993663287883006$& $ a_{64} = -6.07804112569088$& $ a_{65} = -2.820029405964353$\\ $
    a_{71} = -16.625817935606782$& $ a_{72} = -49.999620978741511$& $ a_{73} = 22.3661445506308$\\$
    a_{74} = 30.50526767511958$& $ a_{75} = 13.408435545803448$& $ a_{76} = 1.3455911427944685$ \\
    $b_1 = 0.15881394125505754$& $ b_2 = 3.390357323579911e-13$& $ b_3 = 0.4109696726168125$ \\$
    b_4 = -1.6409254928717294E-13$ & $ b_5 = -0.056173857997504642$ & $ b_6 = 0.40542999348169673$  \\
    $b_7 = 0.08096025064376304$ & $c_2 = 0.34288981581855521$&  $c_3 = 0.2943010566234846$ \\$
    c_4 = 0.34288981581855849$& $ c_5 = -0.028489108793472939$& $ c_6 = 0.73129911961092908$\\$
    c_7 = 1.0000000000000007$ & & \\
    \hline
    \end{tabular}
\end{table}

\clearpage

\bibliographystyle{plain}
\bibliography{references}

\end{document}